\RequirePackage{fix-cm}
\documentclass[smallextended,envcountsect]{svjour3}       
\smartqed  
\usepackage{graphicx}

\usepackage{mathptmx}      
\usepackage{amssymb}

\usepackage{amsmath}
\usepackage{hyperref}
\usepackage{bm}
\usepackage[noadjust]{cite}

\smartqed
\spnewtheorem{thm}[theorem]{Theorem}{\bfseries}{\itshape}
\spnewtheorem{cor}[theorem]{Corollary}{\bfseries}{\itshape}
\spnewtheorem{lem}[theorem]{Lemma}{\bfseries}{\itshape}
\spnewtheorem{ntn}{Notations}[section]{\bfseries}{\itshape}
\spnewtheorem{pro}[theorem]{Proposition}{\bfseries}{\itshape}
\spnewtheorem{dfn}[theorem]{Definition}{\bfseries}{\itshape}
\spnewtheorem{as}{Assumption}[section]{\bfseries}{\itshape}
\spnewtheorem{rem}[theorem]{Remark}{\bfseries}{\itshape}
\spnewtheorem{ob}{Observation}[section]{\bfseries}{\itshape}

\numberwithin{equation}{section}

\newcommand{\C}{\mathbb{C}}    
\newcommand{\N}{\mathbb{N}}    
\newcommand{\PL}{\mathbb{P}}   
\newcommand{\R}{\mathbb{R}}    
\newcommand{\Z}{\mathbb{Z}}    

\newcommand{\wh}{\widehat}
\renewcommand{\le}{\leqslant}
\renewcommand{\ge}{\geqslant}
\newcommand{\bs}{\backslash}
\newcommand{\ol}{\overline}
\newcommand{\la}{\langle}
\newcommand{\ra}{\rangle}
\newcommand{\bo}{\mathcal{O}} 

\newcommand{\pp}{\mathsf{p}}

\newcommand{\DG}{{\mathsf{Diag}}}

\newcommand{\FF}{\mathsf{F}}

\newcommand{\er}{\eqref}

\newcommand{\mra}{{\mathring{a}}}
\newcommand{\mrb}{{\mathring{b}}}
\newcommand{\mrc}{{\mathring{c}}}
\newcommand{\mrd}{{\mathring{d}}}

\newcommand{\mrv}{{\mathring{v}}}
\newcommand{\mrw}{{\mathring{w}}}

\newcommand{\mrta}{{\mathring{\tilde{a}}}}
\newcommand{\mrtb}{{\mathring{\tilde{b}}}}

\newcommand{\mrtv}{{\mathring{\tilde{v}}}}

\newcommand{\mrphi}{\mathring{\phi}}
\newcommand{\mrtvgu}{\mathring{\tilde{\vgu}}}
\newcommand{\bpphi}{\breve{\phi}}
\newcommand{\mrtphi}{\mathring{\tilde{\phi}}}
\newcommand{\mrvgu}{\mathring{\vgu}}

\newcommand{\mrgam}{\mathring{\gamma}}

\newcommand{\tpsi}{\tilde{\psi}}
\newcommand{\tphi}{\tilde{\phi}}

\newcommand{\tvgu}{\tilde{\vgu}}

\newcommand{\bpa}{\breve{a}}
\newcommand{\bpb}{\breve{b}}

\newcommand{\bpta}{\breve{\tilde{a}}}
\newcommand{\bptb}{\breve{\tilde{b}}}

\newcommand{\bpvgu}{\breve{\vgu}}
\newcommand{\bptvgu}{\breve{\tilde{\vgu}}}
\newcommand{\bptphi}{\breve{\tilde{\phi}}}

\newcommand{\htb}{\widehat{\tilde{b}}}

\newcommand{\httheta}{\widehat{\tilde{\theta}}}

\newcommand{\bp}{ \begin{proof} }
	\newcommand{\ep}{\hfill \end{proof} }
\newcommand{\be}{ \begin{equation} }
\newcommand{\ee}{ \end{equation} }
\newcommand{\tp}{\mathsf{T}}
\newcommand{\dm}{\mathsf{M}} 
\newcommand{\dn}{\mathsf{N}} 
\newcommand{\vgu}{\upsilon} 
\newcommand{\Vgu}{\Upsilon} 

\newcommand{\sr}{\operatorname{sr}}  
\newcommand{\vmo}{\operatorname{vm}}
\newcommand{\bvmo}{\operatorname{bvm}}

\newcommand{\td}{\pmb{\delta}}  

\newcommand{\sd}{\mathcal{S}}  
\newcommand{\tz}{\mathcal{T}}  
\newcommand{\cM}{\mathcal{M}}

\newcommand{\cN}{\mathcal{N}}
\newcommand{\bpo}{\operatorname{bo}} 

\newcommand{\om}[1]{\omega_{#1}}
\newcommand{\ga}[1]{\gamma_{#1}}
\newcommand{\ka}[1]{\mathring{\gamma}_{#1}}

\newcommand{\dN}{\mathbb{N}^d}

\newcommand{\dR}{\mathbb{R}^d}

\newcommand{\dZ}{\mathbb{Z}^d}

\newcommand{\dLp}[1]{L_{#1}(\mathbb{R}^d)}

\newcommand{\dLrs}[3]{(L_{#1}(\mathbb{R}^d))^{#2\times #3}}

\newcommand{\vertiii}[1]{{\left\vert\kern-0.25ex\left\vert\kern-0.25ex\left\vert #1 
		\right\vert\kern-0.25ex\right\vert\kern-0.25ex\right\vert}}

\newcommand{\dlp}[1]{l_{#1}(\mathbb{Z}^d)}

\newcommand{\dlrs}[3]{(l_{#1}(\mathbb{Z}^d))^{#2\times #3}}
\newcommand{\dsq}{l(\mathbb{Z}^d)} 

\begin{document}
	
	\title{Compactly supported multivariate dual multiframelets with high vanishing moments and high balancing orders
	}
	
	
	\author{Ran Lu }
	
	
	\institute{Ran Lu \at
		Department of mathematical and statistical sciences, University of Alberta, Edmonton, Canada\\
		\email{rlu3@ualberta.ca} }
	
	\date{Received: date / Accepted: date}

	\maketitle
	
	\begin{abstract}
		Comparing with univariate framelets, the main challenge involved in studying multivariate framelets is that we have to deal with the highly non-trivial problem of factorizing multivariate polynomial matrices. As a consequence, multivariate framelets are much less studied than univariate framelets in the literature. Among existing works on multivariate framelets, multivariate multiframelets are much less considered comparing with the extensively studied scalar framelets. Hence multiframelets are far from being well understood. In this paper, we focus on multivariate dual multiframelets (or dual vector framelets) obtained through the popular oblique extension principle (OEP), which are called OEP-based dual multiframelets. We will show that from any given pair of compactly supported refinable vector functions, one can always construct an OEP-based dual multiframelet, such that its generators have the highest possible order of vanishing moments. Moreover, the associated discrete framelet transform is compact and balanced.
	\end{abstract}
	
	\section{Introduction}
	\label{sec:intro}
Dual framelets derived from refinable vector functions are of interest in applications such as image process and numerical algorithms. The added redundancy in framelet systems enhances their performance over biorthogonal wavelets in practice. For literatures studying framelets/wavelets and their applications, see e.g. \cite{cpss13,cpss15,cs08,ch01,chs02,dh04,dhrs03,dhacha,dh18pp,ds13,Ehler07,eh08,fjs16,goodman94,han97,han03-0,han03,han09,han10,hjsz18,hl19pp,hl20pp,hm03,hm05,jj02,js15,kps16,lj06,lv98,mo06,sz16} and references therein. Dual framelets are usually constructed from \emph{refinable vector functions} via a popular method which is called the \emph{oblique extension principle} (OEP), and such framelets are called \emph{OEP-based framelets}. In this paper, we concentrate on \emph{compactly supported OEP-based dual framelets}. There are three key features which are desired for a compactly supported OEP-based dual framelet in applications: (1) the sparseness of the framelet expansion, which is linked to the vanishing moments of framelet generators; (2) the compactness of the underlying discrete framelet transform, that is, whether or not the transform can be implemented by convolution using finitely supported filters only; (3) the sparseness of the underlying discrete framelet transform, which is closely related to the balancing property of the transform. Quite often, one has to sacrifice (2) to achieve (1), and (3) seems to be too much to expect in most cases. Our goal is to investigate whether or not an OEP-based dual framelet can achieve (1)-(3) simultaneously.

\subsection{Background} 

To better explain our motivations, let us recall some basic concepts. Throughout this paper, $\dm$ is a $d\times d$ dilation matrix, i.e., $\dm\in\Z^{d\times d}$ and its eigenvalues are all greater than one in modulus. For simplicity, let 
$$d_{\dm}:=|\det(\dm)|.$$ 
Denote $ (\dLp{2})^{r\times s}$ the linear space of all $r\times s$ matrices of square integrable functions in $\dLp{2}$. For simplicity, $(\dLp{2})^r:=(\dLp{2})^{r\times 1}$. We introduce the following notion of inner product:
$$\la f,g\ra:=\int_{\dR}f(x)\ol{g(x)}^{\tp}dx,\qquad \forall f\in\dLrs{2}{r}{s},\quad g\in\dLrs{2}{t}{s}.$$
Let $\mrphi,\mrtphi\in (\dLp{2})^r$, $\psi, \tpsi\in (\dLp{2})^s$. We say that $\{\mrphi;\psi\}$ is \emph{an $\dm$-framelet} in $\dLp{2}$ if
there exist positive constants $C_1$ and $C_2$ such that
%
$$C_1\|f\|_{\dLp{2}}^2\le
\sum_{k\in \dZ} |\la f, \mrphi(\cdot-k)\ra|^2+\sum_{j=0}^\infty \sum_{k\in \dZ}
|\la f, \psi_{\dm^j;k}\ra|^2\le C_2\|f\|_{\dLp{2}}^2, \quad  f\in \dLp{2},$$
where $\psi_{\dm^j;k}:=d_{\dm}^{j/2}\psi(\dm^j\cdot-k)$ and $|\la f, \psi_{\dm^j;k}\ra|^2:=\|\la f, \psi_{\dm^j;k}\ra\|^2_{l_2}$. $(\{\mrphi; \psi\},\{\mrtphi;\tpsi\})$ is called \emph{a dual $\dm$-framelet} in $\dLp{2}$ if both $\{\mrphi; \psi\}$ and $\{\mrtphi; \tpsi\}$ are $\dm$-framelets in $\dLp{2}$ and satisfy
\be\label{qtf:expr}
f=\sum_{k\in \dZ} \la f, \mrphi(\cdot-k)\ra \mrtphi(\cdot-k)+
\sum_{j=0}^\infty \sum_{k\in \dZ}
\la f, \psi_{\dm^j;k}\ra \tpsi_{\dm^j;k}, \qquad \forall f\in \dLp{2},
\ee
with the above series converging unconditionally in $\dLp{2}$. $(\{\mrphi; \psi\},\{\mrtphi;\tpsi\})$ is called a dual multiframelet if the multiplicity $r>1$, and is called a scalar framelet if $r=1$. Unless specified, we shall use the term framelet to refer both.

For a dual $\dm$-framelet  $(\{\mrphi;\psi\}, \{\tilde{\mrphi}; \tilde{\psi}\})$, the sparseness of the frame expansion \eqref{qtf:expr} is closely related to the vanishing moments on the framelet generators $\psi$ and $\tpsi$. We say that $\psi$ has $m$ vanishing moments if
$$\int_{\dR}\pp(x)\psi(x)dx=0,\qquad \forall \pp\in\PL_{m-1},$$
where $\PL_{m-1}$ is the space of all $d$-variate polynomials of degree at most $m-1$. Note that $\psi$ has $m$ vanishing moments if and only if
$$\wh{\psi}(\xi)=\bo(\|\xi\|^m),\qquad \xi\to 0,$$
where $f(\xi)=g(\xi)+\bo(\|\xi\|^m)$ as $\xi\to 0$ means $\partial^\mu f(0)=\partial^\mu g(0)$ for all $\mu\in\dN_{0;m}$ with $|\mu|:=\mu_1+\cdots+\mu_d<m$.
We define $\vmo(\psi):=m$ with $m$ being the largest such integer. It is well known in approximation theory (see e.g. \cite[Proposition 5.5.2]{hanbook}) that if $\vmo(\psi)=m$ and $\vmo(\tpsi)=\tilde{m}$, then we necessarily have 
\be \label{qi:order}
\sum_{k\in \dZ} \la \pp, \mrphi(\cdot-k)\ra \mrtphi(\cdot-k):=\sum_{\ell=1}^r \sum_{k\in \dZ} \la \pp, \mrphi_\ell(\cdot-k)\ra \mrtphi_\ell(\cdot-k)=\pp,\qquad \forall\, \pp\in \PL_{m-1}.
\ee
which plays a crucial role in approximation theory and numerical analysis for the convergence rate of the associated approximation/numerical scheme. Moreover, we have
$$\ol{\wh{\mrphi}(\xi)}^\tp \wh{\mrtphi}(\xi+2\pi k)=\bo(\|\xi\|^{m}), \qquad\ol{\wh{\mrphi}(\xi+2\pi k)}^\tp \wh{\mrtphi}(\xi)=\bo(\|\xi\|^{\tilde{m}}),\quad k\in\dZ\bs\{0\},$$
and
$$\ol{\wh{\mrphi}(\xi)}^\tp  \wh{\mrtphi}(\xi)=1+\bo(\|\xi\|^{\tilde{m}+m}),$$
as $\xi\to 0$.

A popular method called \emph{oblique extension principle (OEP)} has been introduced in the literature, which allows us to construct dual framelets with all generators having sufficiently high vanishing moments from \emph{refinable vector functions} \cite{cpss13,cpss15,cs08,ch01,cj00,dhacha,fjs16,hjsz18,js15,kps16,lj06,sz16}. Denote $\dlrs{0}{r}{s}$ the linear space of all $r\times s$ matrix-valued sequences $u=\{u(k)\}_{k\in \dZ}:\dZ\to \C^{r\times s}$ with finitely many non-zero terms. Any element $u\in\dlrs{0}{r}{s}$ is said to be a \emph{finitely supported (matrix-valued) filter/mask}. For $\phi\in (\dLp{2})^r$, we say that $\phi$ is \emph{an $\dm$-refinable vector function} with \emph{a refinement filter/mask} $a\in \dlrs{0}{r}{r}$ if the following \emph{refinement equation} is satisfied:
%
\be\label{ref}\phi(x)=d_{\dm}\sum_{k\in\dZ}a(k)\phi(\dm x-k),\qquad   x\in\dR.\ee
If $r=1$, then we simply say that $\phi$ is an $\dm$-refinable (scalar) function. For $u\in\dlrs{0}{r}{s}$, define its Fourier series via $\wh{u}(\xi):=\sum_{k\in\dZ}u(k)e^{-ik\cdot \xi}$ for $\xi\in\dR$. The Fourier transform is defined via $\wh{f}(\xi):=\int_{\dR} f(x) e^{-ix\cdot\xi} dx$ for $\xi\in \dR$ for all $f\in \dLp{1}$, and can be naturally extended to $\dLp{2}$ functions and tempered distributions. The refinement equation \eqref{ref} is equivalent to
\be\label{ref:f}\wh{\phi}(\dm^{\tp} \xi)=\wh{a}(\xi)\wh{\phi}(\xi),\qquad \xi\in \dR,\ee
where $\wh{\phi}$ is the $r\times 1$ vector obtained by taking entry-wise Fourier transform on $\phi$. Most known framelets are constructed from refinable vector functions via OEP, and we refer them as OEP-based framelets. There are several versions of OEP which have been introduced in the literature (see \cite{chs02,dhrs03,hanbook,hl20pp}). Here we recall the following version of OEP for compactly supported multivariate multiframelets:
\begin{thm}[\emph{Oblique extension principle (OEP)}]\label{thm:df}
	Let $\dm$ be a $d\times d$ dilation matrix. Let $\theta,\tilde{\theta},a,\tilde{a}\in \dlrs{0}{r}{r}$ and $\phi,\tilde{\phi}\in (\dLp{2})^r$ be compactly supported $\dm$-refinable vector functions with refinement filters $a$ and $\tilde{a}$, respectively.
	For matrix-valued filters $b,\tilde{b}\in \dlrs{0}{s}{r}$,
	define
	\be\label{oep:mr}\wh{\mrphi}(\xi):=\wh{\theta}(\xi)\wh{\phi}(\xi),
	\quad \wh{\psi}(\xi):=\htb(\dm^{-\tp}\xi)\wh{\phi}(\dm^{-\tp}\xi),\ee
	\be\label{oep:mrt}\wh{\tilde{\mrphi}}(\xi):=\wh{\tilde{\theta}}(\xi)\wh{\tilde{\phi}}(\xi),
	\quad \wh{\tilde{\psi}}(\xi):=\wh{\tilde{b}}(\dm^{-\tp}\xi)\wh{\tilde{\phi}}(\dm^{-\tp}\xi).
	\ee
	Then $(\{\mrphi;\psi\}, \{\tilde{\mrphi}; \tilde{\psi}\})$ is a dual $\dm$-framelet in $\dLp{2}$ if the following conditions are satisfied:
	\begin{enumerate}
		\item[(1)] $\ol{\wh{\phi}(0)}^\tp \wh{\Theta}(0)\wh{\tilde{\phi}}(0)=1$ with $\wh{\Theta}(\xi):=\ol{\wh{\theta}(\xi)}^\tp \wh{\tilde{\theta}}(\xi)$;
		
		\item[(2)]  $\wh{\psi}(0)=\wh{\tilde{\psi}}(0)=0$.
		
		\item[(3)] $(\{a;b\},\{\tilde{a};\tilde{b}\})_{\Theta}$ forms an OEP-based dual $\dm$-framelet filter bank, i.e.,
		\be \label{dffb}
		\ol{\wh{{a}}(\xi)}^\tp \wh{\Theta}(\dm^{\tp} \xi){\wh{\tilde{a}}(\xi+2\pi \omega)}+\ol{\wh{{b}}(\xi)}^\tp {\wh{\tilde{b}}(\xi+2\pi \omega)}=\td(\omega)\wh{\Theta}(\xi),\ee
		for all $\xi\in\dR$ and $\omega\in\Omega_{\dm}$, where
		\be \label{delta:seq}
		\td(0):=1 \quad \mbox{and}\quad
		\td(x):=0,\qquad \forall\, x\ne 0
		\ee
		and $\Omega_\dm$ is a particular choice of the representatives of cosets in $[\dm^{-\tp}\dZ]/\dZ$ given by
		\be\label{omega:dm}
		\Omega_{\dm}:=\{\om{1},\dots,\om{d_{\dm}}\}:=(\dm^{-\tp}\dZ)\cap [0,1)^d\quad \mbox{with}\quad \om{1}:=0.
		\ee
	\end{enumerate}
\end{thm}
It is clear that the key step to construct an OEP-based dual framelet is to obtain filters $\theta,\tilde{\theta}\in\dlrs{0}{r}{r}$ and $b,\tilde{b}\in\dlrs{0}{s}{r}$ such that $(\{a;b\},\{\tilde{a};\tilde{b}\})_{\Theta}$ is a dual framelet filter bank which satisfies \eqref{dffb}. For any $u\in\dlrs{0}{s}{r}$, define 
\be \label{Pb}
P_{u;\dm}(\xi):=[\wh{u}(\xi+2\pi\om{1}), \ldots, \wh{u}(\xi+2\pi \om{d_{\dm}})],\qquad\xi\in\dR,
\ee
which is an $s\times (rd_{\dm})$ matrix of $2\pi\dZ$-periodic $d$-variate trigonometric polynomials. It is obvious that \eqref{dffb} is equivalent to
\be \label{spectral}
\ol{P_{b;\dm}(\xi)}^\tp P_{\tilde{b};\dm}(\xi)=
\cM_{a,\tilde{a},\Theta}(\xi),
\ee
where 
\be\label{m:a:ta}
\begin{aligned}\cM_{a,\tilde{a},\Theta}(\xi):=&\DG\left(\wh{\Theta}(\xi+2\pi\om{1}),\ldots,
\wh{\Theta}(\xi+2\pi \om{d_{\dm}})\right)\\
&-\ol{P_{a;\dm}(\xi)}^\tp \wh{\Theta}(\dm^\tp\xi) P_{\tilde{a};\dm}(\xi).\end{aligned}\ee

For an OEP-based dual $\dm$-framelet $(\{\mrphi;\psi\}, \{\tilde{\mrphi}; \tilde{\psi}\})$, The orders of vanishing moments of $\psi$ and $\tpsi$ are closely related to the sum rules of the filters $a$ and $\tilde{a}$ associated to $\phi$ and $\tphi$. We say that a filter $a\in (\dlp{0})^{r\times r}$ has \emph{order $m$ sum rules with respect to $\dm$} with a matching filter $\vgu\in \dlrs{0}{1}{r}$ if $\wh{\vgu}(0)\ne 0$ and
\be \label{sr}
\wh{\vgu}(\dm^{\tp}\xi)\wh{a}(\xi+2\pi\omega)=\td(\omega)\wh{\vgu}(\xi)+
\bo(\|\xi\|^m),\quad \xi\to 0,\quad\forall\, \omega\in\Omega_{\dm}.
\ee
In particular, we define 
$$\sr(a,\dm):=\sup\{m\in\N_0: \text{\eqref{sr} holds for some }v\in\dlrs{0}{1}{r}\}.$$
It can be easily deduced from \eqref{dffb} that $\vmo(\psi)\le \sr(\tilde{a},\dm)$ and $\vmo(\tpsi)\le \sr(a,\dm)$ always hold no matter how we choose $\theta$ and $\tilde{\theta}$. Therefore, we are curious about whether or not one can construct filters $\theta,\tilde{\theta}\in\dlrs{0}{r}{r}$ in a way such that the matrix $\cM_{a,\tilde{a},\Theta}$ admits a factorization as in \eqref{spectral} for some $b,\tilde{b}\in\dlrs{0}{s}{r}$ such that $\wh{b}(\xi)\wh{\phi}(\xi)=\bo(\|\xi\|^{\sr(\tilde{a},\dm)})$ and $\wh{\tilde{b}}(\xi)\wh{\tphi}(\xi)=\bo(\|\xi\|^{\sr(a,\dm)})$ as $\xi\to 0$.

\subsection{The major shortcoming of OEP for scalar framelets}

With OEP, a lot of compactly supported scalar dual framelets with the highest possible vanishing moments have been constructed in the literature, to mention only a few, see \cite{cs10,ch00,chs02,dh04,dhrs03,dhacha,dh18pp,han97,han09,hanbook,hm03,hm05,jqt01,js15,mothesis,rs97,sel01} and many references therein. Though OEP appears perfect for improving the vanishing moments of framelet generators, it has a serious shortcoming. To properly address this issue, we need to briefly recall the discrete framelet transform employing an OEP-based filter bank.\\

By $(\dsq)^{s\times r}$ we denote the linear space of all sequences $v: \dZ \rightarrow \C^{s\times r}$. We call every element $v\in(\dsq)^{s\times r}$ \emph{a matrix-valued filter}. For a filter $a\in \dlrs{0}{r}{r}$, we define the filter $a^\star$ via $\wh{a^\star}(\xi):=\ol{\wh{a}(\xi)}^\tp$, or equivalently,
$a^\star(k):=\ol{a(-k)}^\tp$ for all $k\in \dZ$.
We define the \emph{convolution} of two filters via
\[
[v*u](n):=\sum_{k\in \Z} v(k) u(n-k),\quad n\in \dZ,\quad v\in(\dsq)^{s\times r},\quad u\in\dlrs{0}{r}{t}.
\]
Let $\dm$ be a $d\times d$ dilation matrix, define the \emph{upsampling operator} $\uparrow \dm: (\dsq)^{s\times r}\to(\dsq)^{s\times r}$ as
$$[v\uparrow \dm](k):=\begin{cases}v(\dm^{-1} k), &\text{if }k\in \dZ\cap[\dm^{-1}\dZ],\\
0, &\text{elsewhere},\end{cases},\qquad  \forall k\in\dZ,\quad v\in(\dsq)^{s\times r}.$$

We introduce the following operators acting on matrix-valued sequence spaces:

\begin{itemize} 
	
	\item For $u\in\dlrs{0}{r}{t}$, the \emph{subdivision operator} $\sd_{u,\dm}$ is defined via
	$$\sd_{u,\dm} v=|\det(\dm)|^{\frac{1}{2}}[v\uparrow\dm]*u=|\det(\dm)|^{\frac{1}{2}}\sum_{k\in\dZ}v(k)u(\cdot-\dm k),$$
	for all $v\in(\dsq)^{s\times r}$.
	
	\item For $u\in\dlrs{0}{t}{r}$, the \emph{transition operator} $\tz_{u,\dm}$ is defined via
	$$\tz_{u,\dm} v=|\det(\dm)|^{\frac{1}{2}}[v*u^\star]\downarrow\dm =|\det(\dm)|^{\frac{1}{2}}\sum_{k\in\dZ}v(k)\ol{u(k-\dm\cdot)}^{\tp},$$
	for all $v\in(\dsq)^{s\times r}$.
\end{itemize}

Let $\theta,\tilde{\theta},a,\tilde{a}\in\dlrs{0}{r}{r}$ and $b,\tilde{b}\in\dlrs{0}{s}{r}$ be finitely supported filters. For any $J\in\N$ and any input data $v_0\in(\dsq)^{1\times r}$, the $J$-level discrete framelet transform employing the filter bank $(\{a;b\},\{\tilde{a};\tilde{b}\})_{\Theta}$ where $\Theta:=\theta^{\star}*\tilde{\theta}$ is implemented as follows:
\begin{enumerate}
	\item[(S1)] \emph{Decomposition/Analysis}: Recursively compute $v_j,w_j$ for $j=1,\dots,s$ via
	\be\label{dft:anal}v_j:=\tz_{a,\dm}v_{j-1},\qquad w_j:=\tz_{b,\dm}v_{j-1}.\ee
	
	\item[(S2)] \emph{Reconstruction/Synthesis}: Define $\tilde{v}_J:=v_J*\Theta$. Recursively compute $\tilde{v}_{j-1}$ for $j=J,\dots,1$ via
	\be\label{dft:syn}\tilde{v}_{j-1}:=\sd_{\tilde{a},\dm}\mrv_j+\sd_{\tilde{b},\dm}w_j.\ee
	
	\item[(S3)] \emph{Deconvolution}: Recover $\breve{v}_0$ from $\tilde{v}_0$ through $\breve{v}_0*\Theta=\tilde{v}_0.$
\end{enumerate}

We call $\{a;b\}$ the \emph{analysis filter bank} and $\{\tilde{a};\tilde{b}\}$ the \emph{synthesis filter bank}. If any input data $v\in(\dsq)^{1\times r}$ can be exactly retrieved from the above transform, then we say that the $J$-level discrete framelet transform has the \emph{perfect reconstruction property}.

Here comes the major shortcoming of OEP. The deconvolution step (S3) is where the trouble arises. If $(\{a;b\},\{\tilde{a};\tilde{b}\})_{\Theta}$ is an OEP-based dual $\dm$-framelet filter bank satisfying \eqref{dffb}, then the original input data $v_0$ is guaranteed to be a solution of the deconvolution problem $\mathring{v}_0*\Theta=\tilde{v}_0$. However, the deconvolution is inefficient and non-stable, that is, there could be multiple solutions to the deconvolution problem. Thus we cannot expect that the input data can be exactly retrieved by implementing the transform. As observed by \cite[Theorem 2.3]{hl20pp}, a necessary and sufficient condition for a multi-level discrete framelet transform to have the perfect reconstruction property is that $\Theta$ is a \emph{strongly invertible} filter.

\begin{dfn}Let $\Theta\in\dlrs{0}{r}{r}$ be a finitely supported filter. We say that $\wh{\Theta}$ (or simply $\Theta$) is strongly invertible if there exists $\Theta^{-1}\in\dlrs{0}{r}{r}$ such that $\wh{\Theta^{-1}}=\wh{\Theta}^{-1}$, or equivalently all entries of $\wh{\Theta}^{-1}$ are $2\pi\dZ$-periodic trigonometric polynomials.
\end{dfn}

When $\Theta$ is strongly invertible, the discrete framelet transform is said to be \emph{compact}, i.e., the transform is implemented by convolution/deconvolution with finitely supported filters only. The strong invertibility of $\Theta$ forces both $\theta$ and $\tilde{\theta}$ to be strongly invertible. In this case, we can define finitely supported filters $\mra,\mrta\in\dlrs{0}{r}{r}$ and $\mrb,\mrtb\in\dlrs{0}{s}{r}$ via
\be\label{mr:filter}\wh{\mra}(\xi):=\wh{\theta}(\dm^{\tp}\xi)\wh{a}(\xi)\wh{\theta}(\xi)^{-1},\qquad \wh{\mrta}(\xi):=\httheta(\dm^{\tp}\xi)\wh{\tilde{a}}(\xi)\httheta(\xi)^{-1},\ee

\be\label{mrt:filter}\wh{\mrb}(\xi):=\wh{b}(\xi)\wh{\theta}(\xi)^{-1}\qquad \wh{\mrtb}(\xi):=\wh{\tilde{b}}(\xi)\httheta(\xi)^{-1}.\ee
Moreover, if $(\{\mrphi;\psi\},\{\mrtphi;\tpsi\})$ is a dual $\dm$-framelet associated with the OEP-based dual framelet filter bank $(\{a;b\};\{\tilde{a};\tilde{b}\}_{\Theta}$, then the following refinable relations hold:
\be\label{ref:mr}\wh{\mrphi}(\dm^{\tp}\xi)=\wh{\mra}(\xi)\wh{\mrphi(\xi)},\qquad \wh{\psi}(\dm^{\tp}\xi)=\wh{\mrb}(\xi)\wh{\mrphi(\xi)},\ee

\be\label{ref:mrt} \wh{\mrtphi}(\dm^{\tp}\xi)=\wh{\mrta}(\xi)\wh{\mrtphi(\xi)},\qquad \wh{\tpsi}(\dm^{\tp}\xi)=\wh{\mrtb}(\xi)\wh{\mrtphi(\xi)}.\ee
The underlying discrete framelet transform is now employed with the filter bank $(\{\mra;\mrb\},\{ \mrta;\mrtb\})_{I_r}$ without the non-stable deconvolution step as as follows:
\begin{enumerate}
	\item[(S1')] \emph{Decomposition/Analysis}: Recursively compute the \emph{framelet coefficients} $\mrv_j,\mrw_j$ for $j=1,\dots,s$ via
	$$\mrv_j:=\tz_{\mra,\dm}\mrv_{j-1},\qquad \mrw_j:=\tz_{\mrb,\dm}\mrv_{j-1},$$
	where $\mrv_0:=v_0$ is an input data.
	
	\item[(S2')] \emph{Reconstruction/Synthesis}: Define $\mrtv_J:=\mrv_J$. Recursively compute $\mrtv_{j-1}$ for $j=J,\dots,1$ via
	$$\mrtv_{j-1}:=\sd_{\mrta,\dm}\mrtv_j+\sd_{\mrtb,\dm}\mrw_j.$$
	
\end{enumerate}

For a scalar filter $\Theta$ (i.e., $r=1$), it is strongly invertible if and only if $\wh{\Theta}$ is \emph{a non-zero monomial}, i.e.,  $\wh{\Theta}(\xi)=ce^{-ik\cdot\xi}$ for some $c\in\C\setminus\{0\}$ and $k\in\dZ$.  Thus to have a compact discrete framelet transform in the case $r=1$,  $\wh{\theta}$ and $\wh{\tilde{\theta}}$ must be both monomials. However, we lose the main advantage of OEP of improving the vanishing moments of framelet generators by choosing such filters $\theta$ and $\tilde{\theta}$.

\subsection{Advantages and difficulties with multiframelets}
The previously mentioned shortcoming of OEP motivates us to consider multiframelets, that is, framelets with multiplicity $r>1$. Multiframelets have certain advantages over scalar framelets and have been initially studied in \cite{ghm94,glt93} and references therein. In sharp contrast to the extensively studied OEP-based scalar framelets, constructing multiframelets through OEP is much more difficult and is much less studied. To our best knowledge, we are only aware of \cite{han09,hm03,hl19pp,mothesis} for studying one-dimensional OEP-based multiframelets, and \cite{hl20pp} for investigating OEP-based quasi-tight multiframelets in arbitrary dimensions. \\

Here we briefly explain the difficulties involved in studying multiframelets. We see from Theorem~\ref{thm:df} that the most important step of constructing OEP-based framelets is choosing the appropriate filters $\theta,\tilde{\theta}\in\dlrs{0}{r}{r}$. In many situations, this is not easy. Except for the examples in \cite{han09,hl19pp}, all constructed OEP-based dual framelets with non-trivial $\Theta$ (where $\Theta:=\theta^\star*\tilde{\theta}$) do not have a compact underlying discrete framelet transform, i.e., $\Theta$ is not strongly invertible.\\

On the other hand, the sparsity of a discrete framelet transform is another issue which needs to be worried about when the multiplicity $r>1$. First we look at the scalar case $r=1$. Let $(\{\mrphi;\psi\},\{\mrtphi;\tpsi\})$ be an OEP-based dual $\dm$-framelet obtained through Theorem~\ref{thm:df} with an underlying OEP-based dual $\dm$-framelet filter bank $(\{a;b\},\{\tilde{a};\tilde{b}\})_{\Theta}$. Suppose that $\vmo(\psi)=m$. Then the framelet representation \er{qtf:expr} has sparsity in the sense that the polynomial preservation property \er{qi:order} holds. Moreover, item (1) of Theorem~\ref{thm:df} yields $\wh{\phi}(0)\neq 0$. Thus it follows from $\wh{\psi}:=\wh{b}(\dm^{-\tp}\cdot)\wh{\phi}(\dm^{-\tp}\cdot)$ that $\wh{b}(\xi)=\bo(\|\xi\|^m)$ as $\xi\to 0$. For any polynomial $\pp\in\PL_{m-1}$, using Taylor expansion yields $\pp(x-k)=\sum_{\alpha\in\dN_0}\frac{(-k)^{\alpha}}{\alpha!}\partial^{\alpha}\pp(x)$ for all $x,k\in\dR$. Thus for any finitely supported sequence $u\in\dlp{0}$, we have
$$\begin{aligned}[\pp*u](x)&=\sum_{k\in\dZ}\pp(x-k)u(k)=\sum_{\alpha\in\dN_0}[\partial^{\alpha}\pp](x)\left(\sum_{k\in\dZ}\frac{(-k)^{\alpha}}{\alpha!}u(k)\right)\\
&=\sum_{\alpha\in\dN_0}\frac{(-i)^{|\alpha|}}{\alpha!}[\partial^{\alpha}\pp](x)[\partial^{\alpha}\wh{u}](0),\end{aligned}$$
which is a polynomial whose degree is no bigger than the degree of $\pp$, i.e., $\pp*u\in\PL_{m-1}$. Denote $\PL_{m-1}|_{\dZ}$ the linear space of all $d$-variate polynomial sequences of degree at most $m-1$. We now input a polynomial sequence data $\pp\in\PL_{m-1}|_{\dZ}$ and implement the $J$-level discrete framelet transform with the filter bank $(\{a;b\},\{\tilde{a};\tilde{b}\})_{\Theta}$. Observe that the framelet coefficient $v_1$ (see \er{dft:anal}) satisfies $v_1=\tz_{a,\dm}\pp=|\det(\dm)|^{\frac{1}{2}}[\pp*a^{\star}](\dm\cdot)\in\PL_{m-1}$, and by induction we conclude that $v_j\in\PL_{m-1}|_{\dZ}$ for all $j=1,2,\dots, J$. It follows that the framelet coefficients $w_j$ (see \er{dft:syn}) now satisfy
$$w_j=\tz_{b,\dm}v_{j-1}=|\det(\dm)|^{\frac{1}{2}}[v_{j-1}*b^{\star}](\dm\cdot)=\sum_{\alpha\in\dN_0}\frac{(-i)^{|\alpha|}}{\alpha!}[\partial^{\alpha}v_{j-1}](\cdot)[\partial^{\alpha}\wh{b}](0)=0,$$
for all $j=1,\dots,J$, where the last step follows from $\wh{b}(\xi)=\bo(\|\xi\|^m)$ as $\xi\to 0$ and $v_j\in\PL_{m-1}|_{\dZ}$. Consequently, all framelet coefficients $w_j$ vanish. This means that the sparsity of the framelet expansion \er{qtf:expr} automatically guarantees the sparsity of the underlying multi-level discrete framelet transform. Unfortunately this is in general not the case when $r>1$, simply due to the fact that $\wh{\psi}(\xi)=\wh{b}(\dm^{-\tp}\xi)\wh{\phi}(\dm^{-\tp}\xi)=\bo(\|\xi\|^m)$ does not imply any moment property of $\wh{b}(\xi)$ at $\xi=0$. This issue is known as the \emph{balancing property} of a framelet in the literature (\cite{cj00,cj03,han09,han10,hanbook,lv98,sel00}). See Section~\ref{sec:pre} for a brief review of this topic.

\subsection{Main Results and Paper Structure}

From the previous discussion, for OEP-based dual framelets, it seems impossible to achieve high vanishing moments on framelet generators without sacrificing the desired features of the underlying discrete framelet transform. The first breakthrough to this problem is \cite{han09}, which proves that for $r\ge 2$ and $d=1$, one can always obtain OEP-based dual framelets from arbitrary compactly supported refinable vector functions, such that all framelet generators have the highest possible vanishing moments and the associated discrete framelet transform is compact and balanced. However, the case when $d>1$ is far from being well investigated. We are only aware of \cite{han10} which systematically studies the balancing property from the discrete setting for $d>1$ and \cite{hl20pp} which deals with the problem with the approach of the so-called quasi-tight framelets. In this paper, we will systematically study multivariate OEP-based dual framelets with the three key properties. Our main result is the following theorem.

\begin{theorem}\label{thm:dfrt} Let $\dm$ be a $d\times d$ dilation matrix and $r\ge 2$ be an integer. Let $\phi,\tilde{\phi}\in(\dLp{2})^r$ be compactly supported $\dm$ refinable vector functions associated with refinement masks $a,\tilde{a}\in\dlrs{0}{r}{r}$. Suppose that $\sr(a,\dm)=\tilde{m}$ and $\sr(\tilde{a},\dm)=m$ with matching filters $\vgu,\tvgu\in\dlrs{0}{1}{r}$ respectively such that $\wh{\vgu}(0)\wh{\phi}(0)\ne 0$ and $\wh{\tvgu}(0)\wh{\tphi}(0)\ne 0$. Let $\dn$ be a $d\times d$ integer matrix with $|\det(\dn)|=r$. Then there exist $\theta,\tilde{\theta}\in\dlrs{0}{r}{r}$ and $b,\tilde{b}\in\dlrs{0}{s}{r}$ for some $s\in\N$ such that
	
	\begin{enumerate}
		\item[(1)] $\theta$ and $\tilde{\theta}$ are both strongly invertible.
		
		\item[(2)] Define finitely supported filters $\mra,\mrb,\mrta,\mrtb$ via \er{mr:filter} and \er{mrt:filter}. Then $(\{\mra;\mrb\},\{\mrta;\mrtb\})_{I_r}$ is an OEP-based dual $\dm$-framelet filter bank. Moreover, the discrete framelet transform employing the filter bank $(\{\mra;\mrb\},\{\mrta;\mrtb\})_{I_r}$ is order $m$ $E_{\dn}$-balanced, i.e.,  $\bpo(\{\mra;\mrb\},\dm,\dn)=\sr(\mrta,\dm)=m$ (See definition in Section~\ref{sec:pre}).

		\item[(3)] $(\{\mrphi;\psi\},\{\mrtphi;\tpsi\})$ is a compactly supported dual $\dm$-framelet in $\dLp{2}$ with $\vmo(\psi)=m$ and $\vmo(\tpsi)=\tilde{m}$, where $\mrphi,\psi,\mrtphi,\tpsi$ are vector-valued functions defined as in \er{oep:mr} and \er{oep:mrt}.
	\end{enumerate}
	
\end{theorem}

For $r=1$, we have a similar result which only satisfies item (3), for the following reasons: (1) a filter $\theta\in\dlp{0}$ is strongly invertible if and only if $\theta=c\td(\cdot-k)$ for some $c\in\C$ and $k\in\dZ$, and using such filters loses the advantage of OEP for increasing vanishing moments on framelet generators; (2) the balancing property does not come in to play when the multiplicity $r=1$. Theorem~\ref{thm:dfrt} extends the main result of \cite{han09} for the case $d=1$ to $d>1$, but is not a simple generalization. Several techniques for the case $d=1$ simply do not work when $d>1$. For instance, a $2\pi$-periodic trigonometric polynomial has $m$ vanishing moments if and only if it is divisible by $(1-e^{-i\xi})^m$, which is an important fact for the construction of dual framelets with high vanishing moments when $d=1$. Unfortunately, such factorization is no longer available when $d>1$.  A recently developed normal form of a matrix-valued filter (see \cite{hl20pp}) plays a crucial role in our study of OEP-based dual framelets with high vanishing moments and high balancing order, and we will provide a short review of this topic in Section~\ref{sec:pre}.

The paper is organized as follows. In Section~\ref{sec:pre}, we briefly review the balancing property of a multi-level discrete transform, as well as a recently developed normal form of a matrix-valued filter. These are what we need to prove our main result. In Section~\ref{sec:bdft}, we prove the main result Theorem~\ref{thm:dfrt}. Motivated by the proof of the main theorem, we shall perform structural analysis of compactly supported balanced OEP-based dual multiframelets in Section~\ref{sec:str:bdft}. Finally, a summary of our work and some concluding comments will be given in Section~\ref{sec:sum}.

\section{Preliminary}\label{sec:pre}

In this section, we review some important concepts and results which we need to prove our main result on OEP-based dual multiframelets.

\subsection{The balancing property of a multi-level discrete framelet Transform}

As mentioned in Section~\ref{sec:intro}, one issue with OEP when the multiplicity $r>1$ is the sparseness of the multilevel discrete framelet transform. In many applications, the original data is scalar valued, that is, an input data $v\in\dsq$. Thus to implement a multi-level discrete framelet transform, we need to first vectorize the input data. Let $\dn$ be a $d\times d$ integer matrix with $|\det(\dn)|=r$, and let $\Gamma_{\dn}$ be a particular choice of the representatives of the cosets in $\dZ/[\dn\dZ]$ given by 
\be\label{Gamma:N}\Gamma_{\dn}:=\{\mrgam_1,\dots,\mrgam_r\}=:[\dn[0,1)^d]\cap\dZ,\quad\mbox{with}\quad \mrgam_1:=0.\ee
We define the \emph{standard vectorization operator with respect to $\dn$} via
\be\label{vec:con}E_{\dn}v:=(v(\dn\cdot+\mrgam_1),\dots v(\dn\cdot+\mrgam_r)),\qquad\forall v\in\dsq.\ee
Clearly $E_{\dn}$ is a bijection between $\dsq$ and $(\dsq)^{1\times r}$. The sparsity of a multi-level discrete framelet transform employing an OEP-based dual $\dm$-framelet filter bank $(\{a;b\},\{\tilde{a};\tilde{b}\})_{\Theta}$ is measured by the \emph{$E_{\dn}$-balancing order} of the analysis filter bank $\{a;b\}$, denoted by $\bpo(\{a;b\},\dm,\dn):=m$ where $m$ is the largest integer such that the following two conditions hold:
\begin{enumerate}
	\item[(i)] $\tz_{a,\dm}$ is invariant on $E_{\dn}(\PL_{m-1}|_{\dZ})$, i.e., \be \label{bp:lowpass}
	\tz_{a,\dm} E_{\dn}(\PL_{m-1}|_{\dZ})\subseteq E_{\dn}(\PL_{m-1}|_{\dZ}).\ee
	\item[(ii)] The filter $b$ has $m$ \emph{$E_{\dn}$-balancing vanishing moments},  i.e.,
	\be \label{bp:highpass}
	\tz_{b,\dm} E_{\dn} (\pp)=0, \qquad \forall\pp\in \PL_{m-1}|_{\dZ}.
	\ee
	
\end{enumerate}
If items (i) and (ii) are satisfied, note that the framelet coefficient $w_j=\tz_{b,\dm}\tz_{a,\dm}^{j-1}E_{\dn}(\pp)=0$ for all $\pp\in\PL_{m-1}|_{\dZ}$ and $j=1,\dots,J$. This preserves sparsity at all levels of the multilevel discrete framelet transform. A complete characterization of the balancing order of a filter bank is given by the following result.

\begin{thm}\cite[Proposition~3.1, Theorem 4.1]{han10}\label{thm:bp}
	Let $\dm$ be a $d \times d$ dilation matrix and $r\ge 2$ be a positive integer.
	Let $a\in \dlrs{0}{r}{r}$ and $b\in \dlrs{0}{s}{r}$ for some $s\in\N$. Let $\dn$ be a $d\times d$ integer matrix with $|\det(\dn)|=r$ and $E_{\dn}$ in \eqref{vec:con}.
	Define 
	\be \label{vgu:special}
	\wh{\Vgu_{\dn}}(\xi):=\left(e^{i\dn^{-1}\ka{1}\cdot\xi},\ldots, e^{i \dn^{-1}\ka{r}\cdot\xi}\right),\qquad\xi\in\dR.
	\ee
	Then the following statements hold:
	\begin{enumerate}
		
		\item[(1)] The filter $b$ has order $m$ $E_{\dn}$-balancing vanishing moments satisfying \eqref{bp:highpass} if and only if
		\be \label{cond:bvmo} \wh{\Vgu_{\dn}}(\xi)\ol{\wh{b}(\xi)}^\tp
		=\bo(\|\xi\|^m),\qquad \xi \to 0.
		\ee
		\item[(2)] The filter bank $\{a;b\}$ has $m$ $E_{\dn}$-balancing order if and only if \eqref{cond:bvmo} holds and
		%
		$$\wh{\Vgu_{\dn}}(\xi)\ol{\wh{a}(\xi)}^\tp
		=	 \wh{c}(\xi)\wh{\Vgu_{\dn}}(\dm^{\tp}\xi)
		+\bo(\|\xi\|^m),\quad \xi\to 0,$$
for some $c\in \dlp{0}$ with $\wh{c}(0)\ne 0$.
	
	\end{enumerate}
\end{thm}

Let $a,\tilde{a},\theta,\tilde{\theta}\in\dlrs{0}{r}{r}$ and $b,\tilde{b}\in\dlrs{0}{s}{r}$ such that $(\{a;b\},\{\tilde{a};\tilde{b}\})_{\Theta}$ is an OEP-based dual $\dm$-multiframelet filter bank, where $\Theta=\theta^\star*\tilde{\theta}$. Suppose that $\phi,\tilde{\phi}\in (\dLp{2})^r$ are compactly supported $\dm$-refinable vector functions in $\dLp{2}$ satisfying $\wh{\phi}(\dm^{\tp}\xi)=\wh{a}(\xi)\wh{\phi}(\xi)$
and
$\wh{\tilde{\phi}}(\dm^{\tp}\xi)=\wh{\tilde{a}}(\xi)\wh{\tilde{\phi}}(\xi)$.
Define $\mrphi,\psi,\mrtphi,\tpsi$ as in \er{oep:mr} and \er{oep:mrt}.
If $\ol{\wh{\phi}(0)}^\tp \wh{\Theta}(0)\wh{\tilde{\phi}}(0)=1$ and $\wh{\psi}(0)=\wh{\tilde{\psi}}(0)=0$,
then Theorem~\ref{thm:df} tells us that $(\{\mrphi;\psi\},\{\tilde{\mrphi};\tilde{\psi}\})$ is a dual $\dm$-framelet in $\dLp{2}$. With $m:=\sr(\tilde{a},\dm)$, we observe that $ \vmo(\psi)\le m$, $\bvmo(b,\dm,\dn)\le m$ and $\bpo(\{a;b\},\dm,\dn)\le \bvmo(b,\dm,\dn)$.
If $\bpo(\{a;b\},\dm,\dn)=\bvmo(b,\dm,\dn)
=\vmo(\psi)=m$, then we say that the discrete multiframelet transform (or the dual multiframelet $(\{\mrphi;\psi\},\{\tilde{\mrphi};\tilde{\psi}\})$) is \emph{order $m$ $E_{\dn}$-balanced}.
For $r>1$,  $\bpo(\{a,b\},\dm,\dn)<\vmo(\psi)$ often happens. Hence, having high vanishing moments on framelet generators does not guarantee the balancing property and thus significantly reduces the sparsity of the associated discrete multiframelet transform. How to overcome this shortcoming has been extensively studied in the setting of functions in \cite{cj00,lv98,sel00} and in the setting of discrete framelet transforms in \cite{han09,han10,hanbook}.\\

\subsection{The normal form of a matrix-valued filter}

In this section, we briefly review results on a recently developed normal form of the matrix-valued filter. The matrix-valued filter normal form greatly reduces the difficulty in studying multiframelets and multiwavelets, in a way such that we can mimic the techniques we have for studying scalar framelets and wavelets. Considerable works on this topic have been done. We refer the readers to  \cite{han03,han09,han10,hanbook,hl20pp,hm03} for detailed discussion. The most recent advance on this topic is \cite{hl20pp}, which not only generalizes all previously existing works under much weaker conditions but also provides a strengthened normal form of a matrix-valued filter which greatly benefits our study on balanced multivariate multiframelets.

We first recall the following lemma which is known as \cite[Lemma 2.2]{han10}. This result links different vectors of functions which are smooth at the origin by strongly invertible filters.

\begin{lem}\label{lem:U}[\cite[Lemma 2.2]{han10}]
	Let $\wh{v}=(\wh{v_1},\ldots,\wh{v_r})$ and $\wh{u}=(\wh{u_1},\ldots,\wh{u_r})$ be $1\times r$ vectors of functions which are infinitely differentiable at $0$ with $\wh{v}(0)\ne 0$ and $\wh{u}(0)\ne 0$. If $r\ge 2$, then for any positive integer $n\in \N$, there exists a strongly invertible $U\in\dlrs{0}{r}{r}$ such that
	$\wh{u}(\xi)=\wh{v}(\xi)\wh{U}(\xi)+\bo(\|\xi\|^n)$ as $\xi\to 0$.
\end{lem}

One of the most important results on the normal form of a matrix-valued filter is the following result which has been developed recently,  which is a part of \cite[Theorem and 3.3]{hl20pp}.

\begin{thm}\label{thm:normalform:gen}Let $\wh{v},\wh{\mrv}$ be $1\times r$ vectors and $\wh{\phi},\wh{\mrphi}$ be $r\times 1$ vectors of functions which are infinitely differentiable at $0$. Suppose 
	$$\wh{v}(\xi)\wh{\phi}(\xi)=1+\bo(\|\xi\|^m)\quad\mbox{and}\quad \wh{\mrv}(\xi)\wh{\mrphi}(\xi)=1+\bo(\|\xi\|^m),\quad\xi\to 0.$$
	If $r\ge 2$, then for each $n\in\N$, there exists a strongly invertible filter $U\in\dlrs{0}{r}{r}$ such that
	$$\wh{v}(\xi)\wh{U}(\xi)^{-1}=\wh{\mrv}(\xi)+\bo(\|\xi\|^m)\quad\mbox{and}\quad \wh{U}(\xi)\wh{\phi}(\xi)=\wh{\mrphi}(\xi)+\bo(\|\xi\|^n),\quad\xi\to 0.$$
	
\end{thm}

\bp As this result is important for our study of multiframelets but its proof is long and technical, here we provide a sketch of the proof. 

Note that it suffices to prove the claim for $n\ge m$, from which the case $n<m$ follows immediately. The proof contains the following steps:

\begin{enumerate}
	\item[Step 1.] Choose a strongly invertible $U_1\in\dlrs{0}{r}{r}$ such that
	$$\wh{\bpphi}(\xi):=\wh{U_1}(\xi)\wh{\phi}(\xi)=(1,0,\dots,0)^{\tp}+\bo(\|\xi\|^n),\quad \xi\to 0.$$
	Define $\wh{w}:=(\wh{w_1},\dots,\wh{w_r}):=\wh{\vgu}\wh{U_1}^{-1}$ and choose $u_{\vgu}\in\dlrs{0}{1}{r}$ such that $$\wh{u_{\vgu}}(\xi)=(1,\wh{w_2}(\xi),\dots,\wh{w_r}(\xi))\wh{U_1}(\xi)+\bo(\|\xi\|^n),\quad\xi\to 0.$$ 
	Then one can verify that
	$$\wh{u_{\vgu}}(\xi)=\wh{\vgu}(\xi)+\bo(\|\xi\|^m),\quad \wh{u_{\vgu}}(\xi)\wh{\phi}(\xi)=1+\bo(\|\xi\|^n),\quad \xi\to 0.$$
	Similarly we can find $\bptvgu\in\dlrs{0}{1}{r}$ such that
	$$\wh{u_{\mrvgu}}(\xi)=\wh{\mrvgu}(\xi)+\bo(\|\xi\|^m),\quad \wh{u_{\mrvgu}}(\xi)\wh{\mrphi}(\xi)=1+\bo(\|\xi\|^n),\quad \xi\to 0.$$
	
	\item[Step 2.] Choose strongly invertible filters $U_2,U_3\in\dlrs{0}{r}{r}$ such that
	$$(1,0,\dots,0)=\wh{u_{\mrvgu}}(\xi)\wh{U_2}(\xi)+\bo(\|\xi\|^n),\quad \wh{u_{\vgu}}(\xi)=(1,0,\dots,0)\wh{U_3}(\xi)+\bo(\|\xi\|^n),$$
	as $\xi\to 0$. Define $$\wh{u_{\phi}}:=(\wh{u_{\phi,1}},\wh{u_{\phi,2}},\dots,\wh{u_{\phi,r}})^{\tp}:=\wh{U_3}\wh{\phi},\quad \wh{u_{\mrphi}}:=(\wh{u_{\mrphi,1}},\wh{u_{\mrphi,2}},\dots,\wh{u_{\mrphi,r}})^{\tp}:=\wh{U_2}^{-1}\wh{\mrphi}.$$
	It is easy to verify that
	$$\wh{u_{\phi,1}}(\xi)=\wh{u_{\mrphi,1}}(\xi)+\bo(\|\xi\|^n)=1+\bo(\|\xi\|^n),\quad \xi\to 0.$$

	\item[Step 3.] Choose $g_2,\dots,g_r\in\dlp{0}$ such that
	$$\wh{g_{\ell}}(\xi)=\wh{u_{\mrphi,\ell}}(\xi)-\wh{u_{\phi,\ell}}(\xi)+\bo(\|\xi\|^n),\quad\xi\to 0,\quad \ell=2,\dots,r.$$
	Define a strongly invertible filter $U_4\in\dlrs{0}{r}{r}$ via
	$$\wh{U_4}:= \begin{bmatrix}
	1 &0 &\cdots &0\\
	\wh{g_2} &1 &\cdots &0\\
	\vdots &\vdots &\ddots &\vdots\\
	\wh{g_r} &0 &\cdots &1\end{bmatrix}.$$
	Then $U\in\dlrs{0}{r}{r}$ with $\wh{U}:=\wh{U_2}\wh{U_4}\wh{U_3}$ is the desired filter as required.\qed
	
\end{enumerate}\ep

A special case of Theorem~\ref{thm:normalform:gen} is the following result (\cite[Theorem 1.2]{hl20pp}, cf. \cite[Theorem 5.1]{han10}).

\begin{thm}\label{thm:normalform}Let $\dm$ be a $d\times d$ dilation matrix, and let $m\in\N$ and $ r\ge 2$ be integers.
	Let $\phi$ be an $r\times 1$ vector of
	compactly supported distributions satisfying $\wh{\phi}(\dm^{\tp} \xi)=\wh{a}(\xi)\wh{\phi}(\xi)$ with $\wh{\phi}(0)\ne 0$ for some $a\in\dlrs{0}{r}{r}$.
	Suppose the filter $a$ has order $m$ sum rules with respect to $\dm$ satisfying \eqref{sr} with a matching filter $\vgu\in \dlrs{0}{1}{r}$ such that $\wh{\vgu}(0)\wh{\phi}(0)=1$.
	Then for any positive integer $n\in \N$, there exists a strongly invertible filter $U\in\dlrs{0}{r}{r}$ such that the following statements hold:
	\begin{enumerate}
		\item[(1)] Define
		$\wh{\mathring{\vgu}}:=(\wh{\mathring{\vgu}_1},\ldots,	 \wh{\mathring{\vgu}_r}):=\wh{\vgu}\wh{U}^{-1}$ and
		$\wh{\mathring{\phi}}:=(\wh{\mathring{\phi}_1},\ldots,		 \wh{\mathring{\phi}_r})^\tp:=\wh{U}\wh{\phi}$.
		We have
		\be \label{normalform:phi}		 \wh{\mathring{\phi}_1}(\xi)=1+\bo(\|\xi\|^n)
		\quad \mbox{and}\quad \wh{\mathring{\phi}_\ell}(\xi)=\bo(\|\xi\|^n),\quad \xi\to 0,\quad \ell=2,\ldots,r,
		\ee		 \be\label{normalform:vgu}\wh{\mathring{\vgu}_1}(\xi)=1+\bo(\|\xi\|^m)
		\quad \mbox{and}\quad
		\wh{\mathring{\vgu}_\ell}(\xi)=\bo(\|\xi\|^m),\quad \xi\to 0,\quad \ell=2,\ldots,r.\ee	
		
		\item[(2)] Define $\mra\in\dlrs{0}{r}{r}$ via  $\wh{\mathring{a}}:=\wh{U}(\dm^{\tp}\cdot) \wh{a}\wh{U}^{-1}$. Then $\wh{\mrphi}(\dm^{\tp}\cdot)=\wh{\mra}\wh{\mrphi}$ and the new filter $\mra$ has order $m$ sum rules with respect to $\dm$ with the matching filter $\mathring{\vgu}\in (\dlp{0})^{1\times r}$.

	\end{enumerate}
	
\end{thm}

Let $\mra\in\dlrs{0}{r}{r}$ be a refinement mask associated to an $\dm$-refinable vector function $\mrphi$ satisfying \er{normalform:phi}, and suppose that $\mra$ has $m$ sum rules with a matching filter $\mrvgu\in\dlrs{0}{1}{r}$ satisfying \er{normalform:vgu}. It is not hard to observe that $\mra$ has the following structure:
\be \label{normalform}
\wh{\mra}(\xi)=\left[
\begin{matrix}\wh{\mra_{1,1}}(\xi) & \wh{\mra_{1,2}}(\xi)\\
	\wh{\mra_{2,1}}(\xi) & \wh{\mra_{2,2}}(\xi)\end{matrix}\right],
\ee
where $\wh{\mra_{1,1}},\wh{\mra_{1,2}}, \wh{\mra_{2,1}}$ and $\wh{\mra_{2,2}}$ are $1\times 1$, $1\times (r-1)$, $(r-1)\times 1$ and $(r-1)\times (r-1)$ matrices of $2\pi\dZ$-periodic trigonometric polynomials such that
\begin{align}	 &\wh{\mra_{1,1}}(\xi)=1+\bo(\|\xi\|^n),\quad \wh{\mra_{1,1}}(\xi+2\pi\omega)=\bo(\|\xi\|^m),\quad \xi\to 0,\quad\forall \omega\in\Omega_{\dm}\setminus\{0\},\label{mra:11}\\
&\wh{\mra_{1,2}}(\xi+2\pi\omega)=\bo(\|\xi\|^m),\quad \xi\to 0,\quad\forall \omega\in\Omega_{\dm},\label{mra:12}\\ &\wh{\mra_{2,1}}(\xi)=\bo(\|\xi\|^n),\quad \xi\to 0.\label{mra:21}
\end{align}
Any filter $\mra$ satisfying \er{normalform}, \er{mra:11}, \er{mra:12} and \er{mra:21} is said to take the \emph{ideal $(m,n)$-normal form}.\\

If $d=1$, then the three moment conditions \er{mra:11}, \er{mra:12} and \er{mra:21} further yield
$$\wh{\mra_{1,1}}(\xi)=(1+e^{-i\xi}+\dots+e^{-i(|\dm|-1)\xi})^{m}P_{1,1}(\xi)=1+\bo(|\xi|^n),\qquad\xi\to 0,$$
$$\wh{\mra_{1,2}}(\xi)=(1-e^{-i|\dm|\xi})^mP_{1,2}(\xi),\qquad  \wh{\mra_{2,1}}(\xi)=(1-e^{-i\xi})^nP_{2,1}(\xi),$$
where $P_{1,1}, P_{1,2}$ and $P_{2,1}$ are some $1\times 1, 1\times (r-1)$ and $(r-1)\times 1$ matrices of $2\pi$-periodic trigonometric polynomials. Recall that a $2\pi$-periodic trigonometric polynomial $\wh{u}$ satisfies $\wh{u}(\xi)=\bo(|\xi|^m)$ as $\xi\to 0$ if and only if $(1-e^{-i\xi})^m$ divides $\wh{u}$. This is the crucial property to construct univariate dual framelets with high vanishing moments. Unfortunately for $d\ge 2$, there are no corresponding factors for $(1+e^{-i\xi}+\dots+e^{-i(|\dm|-1)\xi})^m$ and $(1-e^{-i\xi})^m$. This means the factorization technique that we have to construct dual framelets with high vanishing moments for the case $d=1$ is no longer available, which illustrates that the investigation is more difficult for $d>1$.

\section{Proof of Theorem~\ref{thm:dfrt}}
\label{sec:bdft}
The goal of this section is to prove the main result Theorem~\ref{thm:df}. To do this, we first need to introduce several notations. For any $k\in\dZ$, the \emph{backward difference operator} $\nabla_k$ is defined via
$$\nabla_ku(n):=u(n)-u(n-k),\qquad \forall n\in\dZ,\quad u\in (\dsq)^{t\times r}.$$
For any multi-index $\alpha:=(\alpha_1,\dots,\alpha_d)\in\dN_0$, we define
$$\nabla^\alpha:=\nabla^{\alpha_1}_{e_1}\nabla^{\alpha_2}_{e_2}\dots\nabla^{\alpha_d}_{e_d},$$
where $\{e_1,\dots,e_d\}$ is the standard basis for $\dR$. Observe that
$$\wh{\nabla^{\alpha}u}(\xi)=
\wh{\nabla^{\alpha}\td}(\xi)\wh{u}(\xi)=
(1-e^{-i\xi_1})^{\alpha_1}(1-e^{-i\xi_2})^{\alpha_2}
\cdots(1-e^{-i\xi_d})^{\alpha_d}\wh{u}(\xi),$$
for all $\xi=(\xi_1,\dots,\xi_d)^{\tp}\in\dR$ and $u\in\dlrs{0}{t}{r}$. 

For $d=1$, recall that a $2\pi$-periodic trignometric polynomial $\wh{c}$ satisfies $\wh{c}(\xi)=\bo(|\xi|^m)$ as $\xi\to 0$ if and only if $\wh{c}$ is divisible by $(1-e^{-i\xi})^m$. Though such a factorization is not available when $d>1$ and there is no factor which plays the role of $(1-e^{-i\xi})^m$ as in the univariate case, the following lemma tells us exactly how one can characterize the moments at zero of a multivariate trigonometric polynomial.

\begin{lem} (\cite[Lemma 5]{dhacha})\label{diff:vm}Let $m\in\N$ and $\wh{c}$ be a $2\pi\dZ$-periodic $d$-variate trigonometric polynomial. Then $\wh{c}(\xi)=\bo(\|\xi\|^m)$ as $\xi\to 0$ if and only if
	$$\wh{c}(\xi)=\sum_{\alpha\in\dN_{0;m}}\wh{\nabla^\alpha\td}(\xi)\wh{c_{\alpha}}(\xi)$$
	for some $c_\alpha\in\dlp{0}$ for all $\alpha\in\dN_{0;m}$, where
	$$\dN_{0;m}:=\{\alpha\in\dN_0:|\alpha|:=\alpha_1+\dots+\alpha_d=m\}.$$
\end{lem}

Next, we introduce the notion of the so-called \emph{coset sequences}. Let $\dm$ be an \emph{invertible integer matrix} and let $\gamma\in\dZ$. For any matrix-valued sequence $u\in(\dsq)^{t\times r}$, we define the \emph{$\gamma$-coset sequence} of $u$ with respect to $\dm$ via
$$u^{[\gamma;\dm]}(k)=u(\gamma+\dm k),\quad k\in\dZ.$$
For $u\in \dlrs{0}{t}{r}$, it is easy to see that
\be\label{f:coset}\wh{u}(\xi)=\sum_{\gamma\in \Gamma_\dm}
\wh{u^{[\gamma;\dm]}}(\dm^{\tp}\xi)e^{-i\gamma\cdot\xi},\qquad\xi\in\dR,\ee
where $\Gamma_\dm$ is a complete set of canonical representatives of the quotient group $\dZ/[\dm\dZ]$, with
\be\label{enum:gamma}
\Gamma_{\dm}:=\{\ga{1},\dots,\ga{d_{\dm}}\}=:(\dm[0,1)^d)\cap\dZ
\quad \mbox{with}\quad \ga{1}:=0.
\ee
Define $\Omega_{\dm}$ via \eqref{omega:dm}. For any filter $u\in\dlrs{0}{r}{r}$ and $\omega\in\Omega_{\dm}$, we introduce the following matrices of trigonometric polynomials associated with $u$ and $\omega$:

\begin{itemize}
	\item Define the $(rd_{\dm} )\times (rd_{\dm})$ block matrix $D_{u,\omega;\dm}(\xi)$, whose $(l,k)$-th $r\times r$ blocks are given by
	\be\label{Duni}(D_{u,\omega;\dm}(\xi))_{l,k}:=
	\begin{cases}\wh{u}(\xi+2\pi\omega), &\text{if }\om{l}+\omega-\om{k}\in\dZ\\
		0, &\text{otherwise}.\end{cases}
	\ee
	
	\item Define the $(rd_{\dm} )\times (rd_{\dm})$ block matrix $E_{u,\omega;\dm}(\xi)$, whose $(l,k)$-th $r\times r$ blocks are given by
	\be\label{Euni}(E_{u,\omega;\dm}(\xi))_{l,k}:=
	\wh{u^{[\ga{k}-\ga{l};\dm]}}
	(\xi)e^{-i\ga{k}\cdot(2\pi\omega)}.
	\ee
	
	\item  Define the $r\times (r d_{\dm})$ matrix $Q_{u;\dm}(\xi)$ via
	\be \label{Qu}
	Q_{u;\dm}(\xi):=\big[\wh{u^{[\ga{1};\dm]}}(\xi),\wh{u^{[\ga{2};\dm]}}(\xi),\dots,\wh{u^{[\ga{d_{\dm}};\dm]}}(\xi)\big].
	\ee
	
\end{itemize}

From \cite[Lemma 7]{dhacha}, it is not hard to deduce that
\be\label{DEF}
\FF_{r;\dm}(\xi)D_{u,\omega;\dm}(\xi)
\ol{\FF_{r;\dm}(\xi)}^{\tp}=d_{\dm} E_{u,\omega;\dm}(\dm^{\tp}\xi),\qquad \xi\in\dR, \omega\in \Omega_{\dm},
\ee	
where $\FF_{r;\dm}(\xi)$ is the following $(rd_{\dm})\times (rd_{\dm})$ matrix:
\be\label{Fourier}
\FF_{r;\dm}(\xi):=\left(e^{-i\ga{l}\cdot(\xi+2\pi\om{k})}I_r\right)_{1\leq l,k\leq d_\dm}.
\ee	
Thus we further deduce that
\be\label{DEF:2}P_{u;\dm}(\xi)=Q_{u;\dm}(\dm^\tp \xi) \FF_{r;\dm}(\xi),\ee
where $P_{u;\dm}(\xi):=\big[\wh{u}(\xi+2\pi\om{1}),
\wh{u}(\xi+2\pi\om{2}),\dots,\wh{u}(\xi+2\pi\om{d_{\dm}})\big]$ as
in \eqref{Pb}. 

Now let $\theta,\tilde{\theta},a,\tilde{a}\in\dlrs{0}{r}{r}$ and $b,\tilde{b}\in\dlrs{0}{s}{r}$ be finitely supported filters. Recall that $(\{a;b\},\{\tilde{a};\tilde{b}\})_{\Theta}$ (where $\Theta:=\theta^\star *\tilde{\theta}$) is a dual $\dm$-framelet filter bank if and only if \er{spectral} holds. Using \er{DEF:2} and $\FF_{r;\dm}\ol{\FF_{r;\dm}}^{\tp}=d_{\dm}I_{d_{\dm}r}$, it is straight forward to see that \er{spectral} is equivalent to
\be\label{dffb:coset:0}
\cN_{a,\tilde{a},\Theta}(\xi)=\ol{Q_{b;\dm}(\xi)}^\tp Q_{\tilde{b};\dm}(\xi),\ee
with
\be\label{dffb:coset}\cN_{a,\tilde{a},\Theta}(\xi):=
d_{\dm}^{-1}E_{\Theta,0;\dm}(\xi)-
\ol{Q_{a;\dm}(\xi)}^\tp \wh{\Theta}(\xi) Q_{\tilde{a};\dm}(\xi).
\ee
Therefore, constructing a dual framelet filter bank is equivalent to obtaining a matrix factorization as in \er{dffb:coset}. When the refinement masks $a$ and $\tilde{a}$ are given, all we have to do is to choose some suitable $\theta$ and $\tilde{\theta}$, and then factorize $\cN_{a,\tilde{a},\Theta}$ as in \er{dffb:coset}. Noting that the matrices $Q_{b;\dm}$ and $Q_{\tilde{b};\dm}$ give us all coset sequences of $b$ and $\tilde{b}$, we can finally reconstruct $b$ and $\tilde{b}$ via \er{f:coset}. It is worth mentioning that the approach of passing to coset sequences often appears in the literature of filter bank construction.\\

Before we prove Theorem~\ref{thm:dfrt}, we need some supporting results. The following result is a special case of \cite[Proposition~3.2]{han03}, which links a refinable vector function $\phi$ with the matching filter $\vgu$ for the associated matrix-valued filter of $\phi$. Here we provide a self-contained simple proof.

\begin{lem}\label{lem:vguphi}
	Let $\dm$ be a $d\times d$ dilation matrix and $a\in \dlrs{0}{r}{r}$.
	Let $\phi$ be an $r\times 1$ vector of
	compactly supported distributions satisfying $\wh{\phi}(\dm^{\tp} \xi)=\wh{a}(\xi)\wh{\phi}(\xi)$ with $\wh{\phi}(0)\ne 0$.
	If $a$ has order $m$ sum rules with respect to $\dm$ satisfying \eqref{sr} with a matching filter $\vgu\in \dlrs{0}{1}{r}$ and $\wh{\vgu}(0)\wh{\phi}(0)=1$, then
	\be \label{vguphi=1:m}
	\wh{\vgu}(\xi)\wh{\phi}(\xi)=1+\bo(\|\xi\|^m),\quad
	\xi\to 0.
	\ee
\end{lem}

\bp  By our assumption on $a$, using $\wh{\vgu}(\dm^{\tp}\xi ) \wh{a}(\xi)=\wh{\vgu}(\xi)+\bo(\|\xi\|^m)$ as $\xi\to 0$ and $\wh{\phi}(\dm^{\tp}\xi)=\wh{a}(\xi)\wh{\phi}(\xi)$, we deduce that
\be \label{rel:vguphi}
\wh{\vgu}(\dm^{\tp} \xi)\wh{\phi}(\dm^{\tp}\xi)
=\wh{\vgu}(\dm^{\tp} \xi) \wh{a}(\xi)\wh{\phi}(\xi)=
\wh{\vgu}(\xi)\wh{\phi}(\xi)+\bo(\|\xi\|^m), \quad \xi\to 0.\ee
We now prove that \er{rel:vguphi} yields \eqref{vguphi=1:m} using \cite[Proposition~2.1]{han03}. For a $p\times q$ matrix $A=(a_{kj})_{1\le k\leq p,1\le q}$ and an $s\times t$ matrix $B$, their Kronecker product $A\otimes B$ is the $(ps) \times (qt)$ block matrix given by
$$
A\otimes B=\begin{bmatrix}a_{11}B & \dots & a_{1q}B\\
\vdots & \ddots & \vdots\\
a_{p1}B & \dots & a_{pq}B\end{bmatrix}.
$$
For any $n\in\N$, define $\otimes^nA:=A\otimes\dots\otimes A$ with $n$ copies of $A$. Recall that if $A,B,C$ and $E$ are matrices of sizes such that one can perform the matrix products $AC$ and $BE$, then we have $(A\otimes B)(C\otimes E)=(AC)\otimes(BE)$. Thus by induction, we have $(\otimes^n (AC))\otimes (BE)=[(\otimes^n A)\otimes B][(\otimes^n C)\otimes E]$.

Define the $1\times d$ vector of differential operators \be\label{diff:op}D:=(\partial_1,\dots,\partial_d),\text{ where }\partial_j:=\frac{\partial}{\partial \xi_j},\quad j=1,\dots, d.\ee 
For simplicity, we define $g(\xi):=\wh{\vgu}(\xi)\wh{\phi}(\xi)$. Direct calculation using the chain rule yields $D\otimes[\wh{g}(\dm^{\tp}\cdot)]=[(D\dm^{\tp})\otimes \wh{g}](\dm^{\tp}\cdot)$. Here $D\dm^{\tp}:=\left(\sum_{j=1}^d\dm_{1j}\partial_j,\dots,\sum_{j=1}^d\dm_{dj}\partial_j\right)$ is a $1\times d$ vector of differential operators where $\dm:=(\dm_{jk})_{1\le j,k\le d}$. By induction, for $j\in \N$, we have
\be\label{mat:kron}[\otimes^{j} D]\otimes[g(\dm^{\tp}\cdot)]=
[(\otimes^{j}(D\dm^{\tp}))\otimes g](\dm^{\tp}\cdot)
=\left([(\otimes^{j}D)\otimes g](\dm^{\tp}\cdot)\right)(\otimes^{j}(\dm^{\tp})).
\ee
It follows from \er{rel:vguphi} and \er{mat:kron} that
\[
\left([(\otimes^{j}D)\otimes g](0)\right)(\otimes^{j}(\dm^{\tp}))
=[(\otimes^{j}D)\otimes g](0),\qquad j=1,\ldots,m-1.
\]
Since all the eigenvalues of $\dm$ are greater than $1$ in modulus, so are the eigenvalues of $\otimes^{j}(\dm^{\tp})$ for every $j\in \N$.
This forces the above linear system to have only the trivial solution $[(\otimes^{j}D)\otimes g](0)=\pmb{0}_{1\times dj}$ for $j=1,\ldots,m-1$.
Hence we conclude that $\partial^\mu g(0)=0$ for all $\mu\in\dN_0$ with $1\le |\mu|\le m-1$. By $g(0)=\wh{\vgu}(0)\wh{\phi}(0)=1$,
we proved $g(\xi)=1+\bo(\|\xi\|^m)$ as $\xi\to 0$, which is just
\eqref{vguphi=1:m}.\qed
\ep

From Theorem~\ref{thm:df}, the most important step for deducing an OEP-based dual multiframelet is choosing suitable filters $\theta,\tilde{\theta}\in\dlrs{0}{r}{r}$ which allow us to perform construction. The following lemma illustrates the existence of $\theta$ and $\tilde{\theta}$ with certain important moment conditions.

\begin{lem}\label{lem:bsr}Let $\dm$ be a $d\times d$ dilation matrix and $r\ge 2$ be an integer. Let $\phi,\tilde{\phi}\in(\dLp{2})^r$ be compactly supported $\dm$ refinable vector functions associated with refinement masks $a,\tilde{a}\in\dlrs{0}{r}{r}$. Suppose that $\sr(a,\dm)=\tilde{m}$ and $\sr(\tilde{a},\dm)=m$ with matching filters $\vgu,\tvgu\in\dlrs{0}{1}{r}$ respectively such that $\wh{\vgu}(0)\wh{\phi}(0)\ne 0$ and $\wh{\tvgu}(0)\wh{\tphi}(0)\ne 0$. Let $\dn$ be a $d\times d$ integer matrix with $|\det(\dn)|=r$, and define $\wh{\Vgu_{\dn}}$ as \er{vgu:special}. Then there exist strongly invertible filters $\theta,\tilde{\theta}\in\dlrs{0}{r}{r}$ such that the following moment conditions hold as $\xi\to 0$:
	\be\label{mr:dual}\wh{\mrvgu}(\xi)=C\ol{\wh{\mrtphi}(\xi)}^{\tp}+\bo(\|\xi\|^{\tilde{m}})=\wh{c}(\xi)\wh{\Vgu_{\dn}}(\xi)+\bo(\|\xi\|^{\tilde{m}}),\ee
	\be\label{mrt:dual}\wh{\mrtvgu}(\xi)=\tilde{C}\ol{\wh{\mrphi}(\xi)}^{\tp}+\bo(\|\xi\|^m)=\wh{d}(\xi)\wh{\Vgu_{\dn}}(\xi)+\bo(\|\xi\|^m),\ee
	\be\label{mrphi:mrtphi}\ol{\wh{\mrphi}(\xi)}^{\tp}\wh{\wh{\mrtphi}}(\xi)=1+\bo(\|\xi\|^{m+\tilde{m}}),\ee	
	for some $c,d\in\dlp{0}$ with $\wh{c}(0)\ne 0$ and $\wh{d}(0)\ne 0$, and some $C,\tilde{C}\in\C\setminus\{0\}$, where $\wh{\mrvgu}:=\wh{\vgu}\wh{\theta}^{-1}, \wh{\mrphi}:=\wh{\theta}\wh{\phi}, \wh{\mrtvgu}:=\wh{\tvgu}\wh{\tilde{\theta}}^{-1}$ and $\wh{\mrtphi}:=\wh{\tilde{\theta}}\wh{\tphi}$.
\end{lem}

\bp By Lemma~\ref{lem:vguphi}, we have
$$\wh{\vgu}(\xi)\wh{\phi}(\xi)=1+\bo(\|\xi\|^{\tilde{m}}),\quad \wh{\tvgu}(\xi)\wh{\tphi}(\xi)=1+\bo(\|\xi\|^{m}),\quad\xi\to 0.$$

Thus by Theorem~\ref{thm:normalform:gen}, there exist strongly invertible filters $\theta,\tilde{\theta}\in\dlrs{0}{r}{r}$ such that 
$$\wh{\mrvgu}(\xi):=\wh{\vgu}(\xi)\wh{\theta}(\xi)^{-1}=r^{-1/2}\wh{\Vgu_{\dn}}(\xi)+\bo(\|\xi\|^{\tilde{m}}),$$ $$\wh{\mrphi}(\xi):=\wh{\theta}(\xi)\wh{\phi}(\xi)=r^{-1/2}\ol{\wh{\Vgu_{\dn}}(\xi)}^{\tp}+\bo(\|\xi\|^n),$$
$$\wh{\mrtvgu}(\xi):=\wh{\tvgu}(\xi)\wh{\tilde{\theta}}(\xi)^{-1}=r^{-1/2}\wh{\Vgu_{\dn}}(\xi)+\bo(\|\xi\|^{m}),$$ $$\wh{\mrtphi}(\xi):=\wh{\tilde{\theta}}(\xi)\wh{\tphi}(\xi)=r^{-1/2}\ol{\wh{\Vgu_{\dn}}(\xi)}^{\tp}+\bo(\|\xi\|^n),$$
as $\xi\to 0$, where $n:=\tilde{m}+m$. This proves \er{mr:dual} and \er{mrt:dual}. Moreover, it is easy to see that \er{mrphi:mrtphi} holds. This completes the proof.\qed\ep

Now we are ready to prove the main result Theorem~\ref{thm:dfrt}.

\textbf{Proof of Theorem~\ref{thm:dfrt}.} By Lemma~\ref{lem:bsr}, there exist strongly invertible filters $\theta,\tilde{\theta}\in\dlrs{0}{r}{r}$ such that \er{mr:dual}, \er{mrt:dual} and \er{mrphi:mrtphi} hold as $\xi\to 0$, where $\wh{\mrvgu}:=\wh{\vgu}\wh{\theta}^{-1}, \wh{\mrphi}:=\wh{\theta}\wh{\phi}, \wh{\mrtvgu}:=\wh{\tvgu}\wh{\tilde{\theta}}^{-1}$ and $\wh{\mrtphi}:=\wh{\tilde{\theta}}\wh{\tphi}$. In particular, we see that item (1) holds.\\

Define $\mra,\mrta\in\dlrs{0}{r}{r}$ as in \er{mr:filter}. We have $\wh{\mrphi}(\dm^{\tp}\cdot)=\wh{\mra}\wh{\mrphi}$, and $\wh{\mrtphi}(\dm^{\tp}\cdot)=\wh{\mrta}\wh{\mrtphi}$. Furthermore, $\mra$ (resp. $\mrta$) has order $\tilde{m}$ (resp. $m$) sum rules with respect to $\dm$ with a matching filter $\mrvgu$ (resp. $\mrtvgu$).\\

Define $n:=\tilde{m}+m$. By Theorem~\ref{thm:normalform:gen}, there exists a strongly invertible $U\in\dlrs{0}{r}{r}$ such that 
$$\wh{\bpphi}(\xi):=\wh{U}(\xi)\wh{\mrphi}(\xi)=(1,0,\dots,0)^{\tp}+\bo(\|\xi\|^n),$$ $$\wh{\bpvgu}(\xi):=\wh{\mrvgu}(\xi)\wh{U}(\xi)^{-1}=(1,0,\dots,0)+\bo(\|\xi\|^{\tilde{m}}),$$
as $\xi\to 0$. Thus by letting $\wh{\bpa}:=\wh{U}(\dm^{\tp}\cdot)\wh{\mra}\wh{U}^{-1}$ , we see that $\bpa$ takes the ideal $(\tilde{m},n)$-normal form, that is, 
$$\wh{\bpa}(\xi)=\begin{bmatrix}\wh{\bpa_{1,1}}(\xi) &\wh{\bpa_{1,2}}(\xi)\\
\wh{\bpa_{2,1}}(\xi) &\wh{\bpa_{2,2}}(\xi)\end{bmatrix},$$
where $\wh{\bpa_{1,1}},\wh{\bpa_{1,2}},\wh{\bpa_{2,1}}$ and $\wh{\bpa_{2,2}}$ are $1\times 1, 1\times (r-1), (r-1)\times 1$ and $(r-1)\times (r-1)$ matrices of $2\pi\dZ$-periodic trigonometric polynomials such that
\begin{align*} &\wh{\bpa_{1,1}}(\xi)=1+\bo(\|\xi\|^n),\quad \wh{\bpa_{1,1}}(\xi+2\pi\omega)=\bo(\|\xi\|^{\tilde{m}}),\quad \xi\to 0,\quad\forall \omega\in\Omega_{\dm}\setminus\{0\},\\
&\wh{\bpa_{1,2}}(\xi+2\pi\omega)=\bo(\|\xi\|^{\tilde{m}}),\quad \xi\to 0,\quad\forall \omega\in\Omega_{\dm},\\ &\wh{\bpa_{2,1}}(\xi)=\bo(\|\xi\|^n),\quad \xi\to 0,
\end{align*}
as $\xi\to 0$, where $\Omega_{\dm}:=\{\omega_1,\dots,\om{d_{\dm}}\}$ is defined as \er{omega:dm}.\\

On the other hand,  we have
$$\wh{\bptvgu}(\xi):=\wh{\mrtvgu}(\xi)\ol{\wh{U}(\xi)}^{\tp}=\ol{\wh{\mrphi}(\xi)}^{\tp}\ol{\wh{U}(\xi)}^{\tp}+\bo(\|\xi\|^m)=(1,0,\dots,0)+\bo(\|\xi\|^m),$$
$$\wh{\bptphi}(\xi):=\ol{\wh{U}(\xi)}^{-\tp}\wh{\mrtphi}(\xi)=\ol{\wh{U}(\xi)}^{-\tp}\ol{\wh{\mrvgu}(\xi)}^{\tp}+\bo(\|\xi\|^{\tilde{m}})=(1,0,\dots,0)^{\tp}+\bo(\|\xi\|^{\tilde{m}}),$$
as $\xi\to 0$. Moreover, the condition \er{mrphi:mrtphi} implies that
$$\wh{\bptphi_1}(\xi)=1+\bo(\|\xi\|^n),\qquad \xi\to 0,$$
where $\bptphi_1$ is the first coordinate of $\bptphi$. Thus by letting $\wh{\bpta}:=\ol{\wh{U}(\dm^{\tp}\cdot)}^{-\tp}\wh{\mrta}\ol{\wh{U}}^{\tp}$, we see that $\wh{\bptphi}(\dm^{\tp}\cdot)=\wh{\bpta}\wh{\bptphi}$ and $\bpta$ has order $m$ sum rules with respect to $\dm$ with a matching filter $\bptvgu$.  Furthermore, we have
$$\wh{\bpta}(\xi)=\begin{bmatrix}\wh{\bpta_{1,1}}(\xi) &\wh{\bpta_{1,2}}(\xi)\\
\wh{\bpta_{2,1}}(\xi) &\wh{\bpta_{2,2}}(\xi)\end{bmatrix},$$ 
where $\wh{\bpta_{1,1}},\wh{\bpta_{1,2}},\wh{\bpta_{2,1}}$ and $\wh{\bpta_{2,2}}$ are $1\times 1, 1\times (r-1), (r-1)\times 1$ and $(r-1)\times (r-1)$ matrices of $2\pi\dZ$-periodic trigonometric polynomials such that
\begin{align*}	 &\wh{\bpta_{1,1}}(\xi)=1+\bo(\|\xi\|^n),\quad \wh{\bpta_{1,1}}(\xi+2\pi\omega)=\bo(\|\xi\|^m),\quad \xi\to 0,\quad\forall \omega\in\Omega_{\dm}\setminus\{0\},\\
&\wh{\bpta_{1,2}}(\xi+2\pi\omega)=\bo(\|\xi\|^m),\quad \xi\to 0,\quad\forall \omega\in\Omega_{\dm},\\ &\wh{\bpta_{2,1}}(\xi)=\bo(\|\xi\|^{\tilde{m}}),\quad \xi\to 0,
\end{align*}
as $\xi\to 0$.\\

For $j=1,\dots,d_{\dm}$, define
$$\wh{A_j}(\xi):=\td(\omega_j)I_r-\ol{\wh{\bpa}(\xi)}^{\tp}\wh{\bpta}(\xi+2\pi\omega_j),$$
where $\td$ is defined as \er{delta:seq}. We have
$$\wh{A_1}(\xi)=I_r-\ol{\wh{\bpa}(\xi)}^{\tp}\wh{\bpta}(\xi)=\begin{bmatrix}\wh{A_{1;1}}(\xi) & \wh{A_{1;2}}(\xi)\\
\wh{A_{1;3}}(\xi) &\wh{A_{1;4}}(\xi)\end{bmatrix},$$
where $\wh{A_{1;1}},\wh{A_{1;2}},\wh{A_{1;3}}$ and $\wh{A_{1;4}}$ are
$1\times 1, 1\times (r-1), (r-1)\times1$ and $(r-1)\times (r-1)$ matrices of $2\pi\dZ$-periodic trigonometric polynomials, satisfying the following moment conditions as $\xi\to 0$:
\begin{align*}
&\wh{A_{1;1}}(\xi)=1-\left(\ol{\wh{\bpa_{1,1}}(\xi)}\wh{\bpta_{1,1}}(\xi)+\ol{\wh{\bpa_{2,1}}(\xi)}^{\tp}\wh{\bpta_{2,1}}(\xi)\right)=\bo(\|\xi\|^n),\\
&\wh{A_{1;2}}(\xi)=-\ol{\wh{\bpa_{1,1}}(\xi)}\wh{\bpta_{1,2}}(\xi)-\ol{\wh{\bpa_{2,1}}(\xi)}^{\tp}\wh{\bpta_{2,2}}(\xi)=\bo(\|\xi\|^{m}),\\
&\wh{A_{1;3}}(\xi)=-\ol{\wh{\bpa_{1,2}}(\xi)}^{\tp}\wh{\bpta_{1,1}}(\xi)-\ol{\wh{\bpa_{2,2}}(\xi)}^{\tp}\wh{\bpta_{2,1}}(\xi)=\bo(\|\xi\|^{\tilde{m}}).\end{align*}
For $j=2,\dots,d_\dm$, we have
$$\wh{A_j}(\xi)=-\ol{\wh{\bpa}(\xi)}^{\tp}\wh{\bpta}(\xi+2\pi\om{j})=\begin{bmatrix}\wh{A_{j;1}}(\xi) & \wh{A_{j;2}}(\xi)\\
\wh{A_{j;3}}(\xi) &\wh{A_{j;4}}(\xi)\end{bmatrix},$$
where $\wh{A_{j;1}},\wh{A_{j;2}},\wh{A_{j;3}}$ and $\wh{A_{j;4}}$ are
$1\times 1, 1\times (r-1), (r-1)\times1$ and $(r-1)\times (r-1)$ matrices of $2\pi\dZ$-periodic trigonometric polynomials for each $j$, satisfying the following moment conditions as $\xi\to 0$:
\begin{align*}
&\wh{A_{j;1}}(\xi)=-\left(\ol{\wh{\bpa_{1,1}}(\xi)}\wh{\bpta_{1,1}}(\xi+2\pi \om{j})+\ol{\wh{\bpa_{2,1}}(\xi)}^{\tp}\wh{\bpta_{2,1}}(\xi+2\pi\om{j})\right)
=\bo(\|\xi\|^m),\\
&\wh{A_{j;1}}(\xi-2\pi \om{j})=-\left(\ol{\wh{\bpa_{1,1}}(\xi-2\pi \om{j})}\wh{\bpta_{1,1}}(\xi)+\ol{\wh{\bpa_{2,1}}(\xi-2\pi \om{j})}^{\tp}\wh{\bpta_{2,1}}(\xi)\right)=\bo(\|\xi\|^{\tilde{m}}),\\
&\wh{A_{j;2}}(\xi)=-\ol{\wh{\bpa_{1,1}}(\xi)}
\wh{\bpta_{1,2}}(\xi+2\pi\om{j})-\ol{\wh{\bpa_{2,1}}(\xi)}^{\tp}\wh{\bpta_{2,2}}(\xi+2\pi\om{j})=\bo(\|\xi\|^m),\\
&\wh{A_{j;3}}(\xi-2\pi\om{j})=-\ol{\wh{\bpa_{1,2}}(\xi-2\pi\om{j})}
\wh{\bpta_{1,1}}(\xi)-\ol{\wh{\bpa_{2,2}}(\xi-2\pi\om{j})}^{\tp}\wh{\bpta_{2,1}}(\xi)= \bo(\|\xi\|^{\tilde{m}}).
\end{align*}

For $\mu\in\dN_0$, define $\Delta_\mu\in \dlrs{0}{r}{r}$ via $\Delta_{\mu}:=\DG(\nabla^{\mu}\td, I_{r-1})$. From what we have done above, we conclude that
\be\label{bpaj:fac:dual}
\wh{A_j}(\xi)=\sum_{\alpha\in\dN_{0;m},\beta\in\dN_{0;{\tilde{m}}}}\ol{\wh{\Delta_\alpha}(\xi)}^{\tp}\wh{A_{j,\alpha,\beta}}(\xi)
\wh{\Delta_\beta}(\xi+2\pi\om{j}),\ee
for some $A_{j,\alpha,\beta}\in\dlrs{0}{r}{r}$ for all $\alpha\in\dN_{0;m},\beta\in\dN_{0;{\tilde{m}}}$ and all $j=1,\dots,d_{\dm}$.\\

Define $\cM_{\bpa,\bpta,I_r}$ as in \er{m:a:ta} with $a,\tilde{a}, \Theta$ being replaced by $\bpa,\bpta,I_r$ respectively, and recall that $D_{\mu,\omega;\dm}$ is defined as \er{Duni} for all $u\in\dlrs{0}{r}{r}$ and $\omega\in\Omega_{\dm}$. Note that $$\cM_{\bpa,\bpta,\tilde{U}}=\sum_{j=1}^{d_{\dm}}D_{A_j,\omega_j;\dm}=\sum_{\alpha\in\dN_{0;m},\beta\in\dN_{0;{\tilde{m}}}}\ol{D_{\Delta_\alpha,0;\dm}}^{\tp}D_{A_{j,\alpha,\beta},\omega_j;\dm}
D_{\Delta_\beta,0;\dm},$$
where the last identity follows from \er{bpaj:fac:dual}. \\

Define $\cN_{\bpa,\bpta,I_r}$ as in \er{dffb:coset} with $a,\tilde{a}$ and $\Theta$ being replaced by $\bpa,\bpta$ and $I_r$ respectively. Recall that $E_{\mu,\omega;\dm}$ is defined as \er{Euni} for all $u\in\dlrs{0}{r}{r}$ and $\omega\in\Omega_{\dm}$, and $\FF_{r;\dm}$ is defined as \er{Fourier}. It follows from \er{DEF:2} and $\FF_{r;\dm}\ol{\FF_{r;\dm}}^{\tp}=d_{\dm}I_{d_{\dm}r}$ that
\be\label{fac:cn}\begin{aligned}&\cN_{\bpa,\bpta,I_r}(\dm^{\tp}\xi)=d_{\dm}^{-2}\FF_{r;\dm}(\xi)\cM_{\bpa,\bpta,I_r}(\xi)\ol{\FF_{r;\dm}(\xi)}^{\tp}\\
	=&d_{\dm}^{-1}\sum_{j=1}^{d_{\dm}}
	\sum_{\alpha\in\dN_{0;m},\beta\in\dN_{0;\tilde{m}}}
	\ol{E_{\Delta_\alpha,0;\dm}(\dm^{\tp}\xi)}^{\tp}
	E_{A_{j,\alpha,\beta},\omega_j;\dm}(\dm^{\tp}\xi)
	E_{\Delta_\beta,0;\dm}(\dm^{\tp}\xi).
\end{aligned}\ee
By letting
$$E_{\alpha,\beta}(\xi):=d_{\dm}^{-1}\sum_{j=1}^{d_{\dm}}E_{A_{j,\alpha,\beta},\omega_j;\dm}(\xi),\qquad \xi\in\dR,\quad\alpha\in\dN_{0;m},\quad\beta\in\dN_{0;\tilde{m}}, $$
we have
\be\label{fac:cn:2}\cN_{\bpa,\bpta,I_r}(\xi)=
\sum_{\alpha\in\dN_{0;m},\beta\in\dN_{0;\tilde{m}}}
\ol{E_{\Delta_\alpha,0;\dm}(\xi)}^{\tp}
E_{\alpha,\beta}(\xi)
E_{\Delta_\beta,0;\dm}(\xi).\ee
For every $\alpha\in\dN_{0;m}$ and $\beta\in\dN_{0;\tilde{m}}$, choose $E_{\alpha,\beta,1}$ and $E_{\alpha,\beta,1}$ which are $d_{\dm}r \times d_{\dm}r$ matrices of $2\pi\dZ$-periodic trigonometric polynomials such that $E_{\alpha,\beta}=\ol{E_{\alpha,\beta,1}}^{\tp}E_{\alpha,\beta,2}$. Define $\bpb_{\alpha,\beta, k},\bptb_{\alpha,\beta,k}\in\dlrs{0}{1}{r}$ for $k=1,\dots, d_{\dm}r$ and all $\alpha\in\dN_{0;m},\beta\in\dN_{0;\tilde{m}}$ via
\begin{align}&\wh{\bpb_{\alpha,\beta}}(\xi):=\begin{bmatrix}\wh{\bpb_{\alpha,\beta,1}}(\xi)\\
\vdots\\
\wh{\bpb_{\alpha,\beta,d_{\dm}r}}(\xi)\end{bmatrix}:=E_{\alpha,\beta,1}(\dm^{\tp}\xi)\FF_{r;\dm}(\xi)\begin{bmatrix}\wh{\Delta_{\alpha}}(\xi)\\
\pmb{0}_{d_{\dm}(r-1)\times r}\end{bmatrix},\label{bab}\\
&\wh{\bptb_{\alpha,\beta}}(\xi):=\begin{bmatrix}\wh{\bptb_{\alpha,\beta,1}}(\xi)\\
\vdots\\
\wh{\bptb_{\alpha,\beta,d_{\dm}r}}(\xi)\end{bmatrix}:=E_{\alpha,\beta,2}(\dm^{\tp}\xi)\FF_{r;\dm}(\xi)\begin{bmatrix}\wh{\Delta_{\beta}}(\xi)\\
\pmb{0}_{d_{\dm}(r-1)\times r}\end{bmatrix},\label{btab}
\end{align}
where $\pmb{0}_{t\times q}$ denotes the $t\times q$ zero matrix. Recall that $P_{u;\dm}(\xi)=[\wh{u}(\xi+2\pi\omega_1),\dots,\wh{u}(\xi+2\pi\omega_{d_{\dm}})]$ as in \er{Pb} for all matrix-valued filter $u$. It is not hard to observe that
\be\label{Pbab}\begin{aligned}P_{\bpb_{\alpha,\beta};\dm}(\xi)&=E_{\alpha,\beta,1}(\dm^{\tp}\xi)\FF_{r;\dm}(\xi)D_{\Delta_{\alpha},0;\dm}(\xi)\\
	&=E_{\alpha,\beta,1}(\dm^{\tp}\xi)E_{\Delta_{\alpha},0;\dm}(\dm^{\tp}\xi)\ol{\FF_{r;\dm}(\xi)}^{\tp},\end{aligned}\ee 
where the last identity follows from \er{DEF} and $\FF_{r;\dm}\ol{\FF_{r;\dm}}^{\tp}=d_{\dm}I_{d_{\dm}r}$. Similarly,
\be\label{Pbtab}P_{\bptb_{\alpha,\beta};\dm}(\xi)=E_{\alpha,\beta,2}(\dm^{\tp}\xi)E_{\Delta_{\beta},0;\dm}(\dm^{\tp}\xi)\ol{\FF_{r;\dm}(\xi)}^{\tp}.\ee
It follows from \er{fac:cn}, \er{fac:cn:2}, \er{Pbab} and \er{Pbtab} that
\be\label{fac:cm}\begin{aligned}\cM_{\bpa,\bpta,I_r}(\xi)&=\ol{\FF_{r;\dm}(\xi)}^{\tp}\cN_{\bpa,\bpta,I_r}(\dm^{\tp}\xi)\FF_{r;\dm}(\xi)\\
&=\sum_{\alpha\in\dN_{0;m},\beta\in\dN_{0;\tilde{m}}}\ol{P_{\bpb_{\alpha,\beta};\dm}(\xi)}^{\tp}P_{\bptb_{\alpha,\beta};\dm}(\xi).
\end{aligned}\ee
Define
\begin{align*}&\{\bpb_{\ell}:\ell=1,\dots,s\}:=\{\bpb_{\alpha,\beta}:\alpha\in\dN_{0;m},\quad\beta\in\dN_{0;\tilde{m}}\},\\
&\{\bptb_{\ell}:\ell=1,\dots,s\}:=\{\bptb_{\alpha,\beta}:\alpha\in\dN_{0;m},\quad\beta\in\dN_{0;\tilde{m}}\},\end{align*}
and let $\bpb:=[\bpb_1^{\tp},\dots,\bpb_s^{\tp}]^{\tp}, \bptb:=[\bptb_1^{\tp},\dots,\bptb_s^{\tp}]^{\tp}.$ We see that \er{fac:cm} becomes
$$\cM_{\bpa,\bpta,I_r}(\xi)=\ol{P_{\bpb;\dm}(\xi)}^{\tp}P_{\bptb;\dm}(\xi),$$
which is equivalent to say that $(\{\bpa;\bpb\},\{\bpta;\bptb\})_{I_r}$ is an OEP-based dual $\dm$-framelet filter bank satisfying
\be\label{oep:bp}\ol{\wh{\bpa}(\xi)}^{\tp}\wh{\bpta}(\xi+2\pi\omega)+\ol{\wh{\bpb}(\xi)}^{\tp}\wh{\bptb}(\xi+2\pi\omega)=\td(\omega)I_r,\qquad \xi\in\dR,\omega\in\Omega_{\dm}.\ee
Now define $\mrb,\mrtb,b,\tilde{b}\in\dlrs{0}{s}{r}$ via
$$\wh{\mrb}:=\wh{\bpb}\wh{U}^{-1},\quad \wh{\mrtb}:=\wh{\bptb}\ol{\wh{U}}^{\tp},\quad \wh{b}:=\wh{\mrb}\wh{\theta}^{-1},\quad \wh{\tilde{b}}:=\wh{\mrtb}\wh{\tilde{\theta}}^{-1}.$$
It follows from \er{oep:bp} that $(\{\mra;\mrb\},\{\mrta;\mrtb\})_{I_r}$ is an OEP-based dual $\dm$-framelet filter bank satisfying
$$\ol{\wh{\mra}(\xi)}^{\tp}\wh{\mrta}(\xi+2\pi\omega)+\ol{\wh{\mrb}(\xi)}^{\tp}\wh{\mrtb}(\xi+2\pi\omega)=\td(\omega)I_r,\qquad \xi\in\dR,\omega\in\Omega_{\dm},$$
and $(\{a;b\},\{\tilde{a};\tilde{b}\})_{\Theta}$ (where $\Theta:=\theta^{\star}*\tilde{\theta}$) is an OEP-based dual $\dm$-framelet filter bank satisfying \er{dffb}. By \er{mrt:dual} and \er{bab}, we have
\be\label{bpo:b}\begin{aligned}&\wh{\Vgu_{\dn}}(\xi)\ol{\wh{\mrb}(\xi)}^{\tp}=\ol{\wh{d}(\xi)^{-1}\wh{\mrphi}(\xi)}^{\tp}\ol{\wh{\mrb}(\xi)}^{\tp}+\bo(\|\xi\|^m)=\ol{\wh{d}(\xi)}^{-1}\ol{\wh{\bpphi}(\xi)}^{\tp}\ol{\wh{\bpb}(\xi)}^{\tp}+\bo(\|\xi\|^m)\\
	=&\ol{\wh{d}(\xi)}^{-1}(1,0,\dots,0)\ol{\wh{\bpb}(\xi)}^{\tp}+\bo(\|\xi\|^m)=\bo(\|\xi\|^m),\quad \xi\to 0,\end{aligned}\ee
where $d\in\dlp{0}$ with $\wh{d}(0)\ne 0$ is the same as in \er{mrt:dual}. 
Similarly, we deduce from \er{mr:dual} and \er{btab} that
\be\label{bpo:b:3}\wh{\Vgu_{\dn}}(\xi)\ol{\wh{\mrtb}(\xi)}^{\tp}=\bo(\|\xi\|^{m}),\qquad\xi\to 0.\ee
On the other hand, it follows immediately from \er{mrt:dual} and the refinement relation $\wh{\mrphi}(\dm^{\tp}\cdot)=\wh{\mra}\wh{\mrphi}$ that
$$\frac{\wh{d}(\xi)}{\wh{d}(\dm^{\tp}\xi)}\wh{\Vgu_{\dn}}(\xi)\ol{\wh{\mra}(\xi)}^{\tp}=\wh{\Vgu_{\dn}}(\dm^{\tp}\xi)+\bo(\|\xi\|^{\tilde{m}}),\qquad \xi\to 0.$$
Hence by Theorem~\ref{thm:bp}, we have $\bpo(\{\mra;\mrb\},\dm,\dn)=m=\sr(\mrta;\dm)$. This proves item (2).\\

Now define vector functions $\psi$ and $\psi$ as in \er{oep:mr} and \er{oep:mrt}. It follows from \er{mr:dual}, \er{mrt:dual}, \er{bpo:b} and \er{bpo:b:3} that $\vmo(\psi)=m$ and $\vmo(\tpsi)=\tilde{m}$. Further note that
$$\ol{\wh{\phi}(0)}^{\tp}\wh{\Theta}(0)\wh{\tphi}(0)=\ol{\wh{\mrphi}(0)}^{\tp}\wh{\mrtphi}(0)=1.$$
It follows from Theorem~\ref{thm:df} that $(\{\mrphi;\psi\},\{\mrtphi;\tpsi\})$ is a dual $\dm$-framelet in $\dLp{2}$. This proves item (3).\qed

Theorem~\ref{thm:dfrt} is valid for the case $r>1$. For the case $r=1$, we have to sacrifice the strong invertibility of $\theta$ and $\tilde{\theta}$ to improve the orders of vanishing moments of the framelet generators. Nevertheless, the matrix decomposition technique in the proof of Theorem~\ref{thm:dfrt} can be applied to deduce the following result for the case $r=1$.

\begin{cor}\label{cor:dft:r:1} Let $\dm$ be a $d\times d$ dilation matrix and let $\phi,\tphi\in \dLp{2}$ be compactly supported refinable functions satisfying $\wh{\phi}(\dm^{\tp} \xi)=\wh{a}(\xi)\wh{\phi}(\xi)$ and $\wh{\tphi}(\dm^{\tp} \xi)=\wh{\tilde{a}}(\xi)\wh{\tphi}(\xi)$, where $a,\tilde{a}\in\dlp{0}$ have order $\tilde{m}$ and $m$ sum rules with respect to $\dm$ with matching filters $\vgu,\tvgu\in \dlp{0}$, respectively. Suppose that $\wh{\vgu}(0)\wh{\phi}(0)=\wh{\tvgu}(0)\wh{\tphi}(0)=1$. Then there exist $b,\tilde{b}\in \dlrs{0}{s}{1}$ and $\theta,\tilde{\theta}\in \dlp{0}$ such that
	
	\begin{enumerate}
		\item  $(\{a;b\},\{\tilde{a};\tilde{b}\})_{\theta^\star*\tilde{\theta}}$ forms an OEP-based dual $\dm$-framelet filter bank.

		\item $(\{\mrphi;\psi\},\{\mrtphi;\tpsi\})$ is a compactly supported dual $\dm$-framelet in $\dLp{2}$, where $\mrphi,\psi,\mrtphi$ and $\tpsi$ are defined as in \er{oep:mr} and \er{oep:mrt}. Moreover, $\vmo(\psi)=m$ and $\vmo(\tpsi)=\tilde{m}$.
	
	\end{enumerate}
\end{cor}

\section{Structural investigation on balanced OEP-based dual framelets}\label{sec:str:bdft}

In this section, we perform structural analysis on OEP-based dual framelets with hight balancing orders. 

The most important step to obtain an OEP-based dual framelet with high balancing orders is finding the suitable filters $\theta,\tilde{\theta}$. From Theorem~\ref{thm:dfrt} and its proof, we have some clue on the choices of such filters. The following theorem states the sufficient conditions for obtaining an OEP-based dual framelet with all desired properties.

\begin{theorem}\label{theta:1}Let $\dm$ be a $d\times d$ dilation matrix and $r\ge 2$ be an integer. Let $\phi,\tilde{\phi}\in(\dLp{2})^r$ be compactly supported $\dm$ refinable vector functions associated with refinement masks $a,\tilde{a}\in\dlrs{0}{r}{r}$. Suppose that $\sr(a,\dm)=\tilde{m}$ and $\sr(\tilde{a},\dm)=m$ with matching filters $\vgu,\tvgu\in\dlrs{0}{1}{r}$ respectively such that $\wh{\vgu}(0)\wh{\phi}(0)\ne 0$ and $\wh{\tvgu}(0)\wh{\tphi}(0)\ne 0$. Let $\dn$ be a $d\times d$ integer matrix with $|\det(\dn)|=r$, and define $\wh{\Vgu_{\dn}}$ as in \er{vgu:special}.

	Let $\theta,\tilde{\theta}\in\dlrs{0}{r}{r}$ be strongly invertible finitely supported filters. Then
	
	\begin{enumerate}
		\item[(i)] the moment conditions \er{mr:dual},\er{mrt:dual} and \er{mrphi:mrtphi} hold as $\xi\to 0$, for some $c,d\in\dlp{0}$ with $\wh{c}(0)\ne 0$ and $\wh{d}(0)\ne 0$, and some $C,\tilde{C}\in\C\setminus\{0\}$, where $\wh{\mrvgu}:=\wh{\vgu}\wh{\theta}^{-1},\wh{\mrphi}:=\wh{\theta}\wh{\phi},\wh{\mrtvgu}:=\wh{\tvgu}\wh{\tilde{\theta}}^{-1}$ and $\wh{\mrtphi}:=\wh{\tilde{\theta}}\wh{\tphi}$,
	\end{enumerate}
	
	implies
	
	\begin{enumerate}
		\item[(ii)] there exist finitely supported filters $b,\tilde{b}\in\dlrs{0}{s}{r}$ such that all claims in Theorem~\ref{thm:dfrt} hold.
	\end{enumerate}	 
	
	Conversely, if in addition assume that
	
	\begin{enumerate}
		\item[(iii)] $1$ is a simple eigenvalue of $\wh{a}(0)$ and $\wh{\mrta}(0)$. Moreover,  
		$$\lambda^{\alpha}I_r-\wh{a}(0),\quad I_r -\lambda^{\beta}\wh{a}(0),\quad I_r-\lambda^{\alpha}\wh{\tilde{a}}(0),\quad \lambda^{\beta}I_r -\wh{\tilde{a}}(0) $$ 
		are invertible matrices for all $\alpha, \beta\in\dN_0$ with $0<|\alpha|<\tilde{m}$ and $0<|\beta|<m$, where $\lambda:=(\lambda_1,\dots,\lambda_d)$ is the vector of the eigenvalues of $\dm$.

		\item[(iv)] $\wh{p}(\dm^{\tp}\xi)\wh{\Vgu_{\dn}}(\dm^{\tp}\xi)\wh{\mrta}(\xi)=\wh{p}(\xi)\wh{\Vgu_{\dn}}(\xi)+\bo(\|\xi\|^{m})$ as $\xi\to 0$ for some $p\in\dlp{0}$ with $\wh{p}(0)\neq 0$, where $\wh{\mrta}:=\wh{\tilde{\theta}}(\dm^{\tp}\cdot)\wh{\tilde{a}}\wh{\tilde{\theta}}^{-1}$ is defined as in \er{mrt:filter}.
		
		\item[(v)] $\wh{q}(\xi)\wh{\mrta}(\xi)\ol{\wh{\Vgu_{\dn}}(\xi)}^{\tp}=\wh{q}(\dm^{\tp}\xi)\ol{\wh{\Vgu_{\dn}}(\dm^{\tp}\xi)}^{\tp}+\bo(\|\xi\|^{\tilde{m}})$ as $\xi\to 0$ for some $q\in\dlp{0}$ with $\wh{q}(0)\neq 0$.

	\end{enumerate}
	Then item (ii) implies (i).	
	
\end{theorem}

\bp  The implication (i) $\Rightarrow$ (ii) follows immediately from the proof of Theorem~\ref{thm:dfrt}.

Now suppose item (ii) holds. Define $\mra,\mrta\in\dlrs{0}{r}{r}$ as in \er{mr:filter} and define
$\mrb,\mrtb\in\dlrs{0}{s}{r}$ as in \er{mrt:filter}. By item (2) of Theorem~\ref{thm:dfrt}, we have
\be\label{dffrt}\ol{\wh{\mra}(\xi)}^{\tp}\wh{\mrta}(\xi)+\ol{\wh{\mrb}(\xi)}^{\tp}\wh{\mrtb}(\xi)=I_r,\ee
and $\bpo(\{\mra;\mrb\},\dm,\dn)=m$. By Theorem~\ref{thm:bp}, we have
\be\label{bp:mr}\wh{\Vgu_{\dn}}(\xi)\ol{\wh{\mrb}(\xi)}^{\tp}=\bo(\|\xi\|^m),\quad \wh{\Vgu_{\dn}}(\xi)\ol{\wh{\mra}(\xi)}^{\tp}=\wh{\mrc}(\xi)\wh{\Vgu_{\dn}}(\dm^{\tp}\xi)+\bo(\|\xi\|^m),\quad \xi\to 0,\ee
for some $\mrc\in\dlp{0}$ with $\wh{\mrc}(0)\neq 0$.

Assume in addition that items (iii) - (v) hold.

By left multiplying $\wh{\Vgu_{\dn}}$ on both sides of \er{dffrt} and using item (iv), we have
$$\wh{\Vgu_{\dn}}(\xi)=\wh{\mrc}(\xi)\wh{\Vgu_{\dn}}(\dm^{\tp}\xi)\wh{\mrta}(\xi)+\bo(\|\xi\|^m)=\wh{\mrc}(\xi)\frac{\wh{p}(\xi)}{\wh{p}(\dm^{\tp}\xi)}\wh{\Vgu_{\dn}}(\xi)+\bo(\|\xi\|^m),\quad\xi\to 0.$$
From the above relation we conclude that $\wh{\mrc}(0)=1$, and thus
\be\label{sr:mra}\wh{\mrd}(\dm^{\tp}\xi)\wh{\Vgu_{\dn}}(\dm^{\tp}\xi)\wh{\mrta}(\xi)=\wh{\mrd}(\xi)\wh{\Vgu_{\dn}}(\xi)+\bo(\|\xi\|^{m}),\quad \xi\to 0,\ee
where $\mrd\in\dlp{0}$ satisfies 
$$\wh{\mrd}(\xi)=\prod_{j=1}^\infty\wh{\mrc}((\dm^{\tp})^{-j}\xi)+\bo(\|\xi\|^m),\qquad\xi\to 0.$$
Moreover, it is easy to see from the second relation in \er{bp:mr} that 
\be\label{ref:mrta}\ol{\wh{\mrd}(\dm^{\tp}\xi)\wh{\Vgu_{\dn}}(\dm^{\tp}\xi)}^{\tp}=\wh{\mra}(\xi)\ol{\wh{\mrd}(\xi)\wh{\Vgu_{\dn}}(\xi)}^{\tp}+\bo(\|\xi\|^m),\quad\xi\to 0.\ee

We now apply the argument in the proof of \cite[Lemma 2.2]{han03} to prove that \er{mr:dual} and \er{mrt:dual} must hold.

Since $\mrta$ has $m$ sum rules with a matching filter $\mrtvgu$ with $\wh{\mrtvgu}:=\wh{\tvgu}\wh{\tilde{\theta}}^{-1}$, we have $\wh{\mrtvgu}(\dm^{\tp}\xi)\wh{\mrta}(\xi)=\wh{\mrtvgu}(\xi)+\bo(\|\xi\|^{m})$ as $\xi\to 0$. This implies that
\be\label{tvgu:mom}\left[(\otimes^jD)\otimes\left(\wh{\mrtvgu}(\dm^{\tp}\cdot)\wh{\mrta}\right)\right](0)=[(\otimes^jD)\otimes\wh{\mrtvgu}](0),\qquad j=1,\dots,m-1,\ee
where $D$ is the vector of differential operators defined as \er{diff:op}. Rearranging \er{tvgu:mom} yields
\be\label{tvgu:mom:1}\left[(\otimes^jD)\otimes\left(\wh{\mrtvgu}(\dm^{\tp}\cdot)\wh{\mrta}(0)-\wh{\mrtvgu}\right)\right](0)=\left[(\otimes^jD)\otimes\left(\wh{\mrtvgu}(\dm^{\tp}\cdot)(\wh{\mrta}(0)-\wh{\mrta})\right)\right](0),\ee
for all $j=1,\dots,m-1.$ By the generalized product rule, we observe that the right hand of \er{tvgu:mom:1} only involves $\partial^\mu\wh{\mrtvgu}(0)$ with $|\mu|<j$. By calculation, we have
\be\label{tvgu:mom:2}\begin{aligned}&\left[(\otimes^jD)\otimes\left(\wh{\mrtvgu}(\dm^{\tp}\cdot)\wh{\mrta}(0)-\wh{\mrtvgu}\right)\right](0)\\
	=&\left([(\otimes^jD)\otimes \wh{\mrtvgu}](0)\right)[(\otimes^j\dm^{\tp})\otimes \wh{\mrta}(0)-I_{d^jr}],\end{aligned}\ee
for all $j\in\N$. Now by the condition in item (iii) on $\mrta$, the matrix $[(\otimes^j\dm^{\tp})\otimes \wh{\mrta}(0)-I_{d^jr}]$ is invertible for $j=1,\dots,m-1$. Moreover, it follows from \er{tvgu:mom:1} and \er{tvgu:mom:2} that up to a multiplicative constant, the values $\partial^\mu\wh{\mrtvgu}(0), |\mu|<m$ are uniquely determined via $\wh{\mrtvgu}(0)\wh{\mrta}(0)=\wh{\mrtvgu}(0)$ and
$$[(\otimes^jD)\otimes \wh{\mrtvgu}](0)=\left(\left[(\otimes^jD)\otimes\left(\wh{\mrtvgu}(\dm^{\tp}\cdot)(\wh{\mrta}(0)-\wh{\mrta})\right)\right](0)\right)[(\otimes^j\dm^{\tp})\otimes \wh{\mrta}(0)-I_{d^jr}]^{-1},$$
for all $j=1,\dots,m-1$.

Next, note that the refinement relation $\wh{\mrphi}(\dm^{\tp}\cdot)=\wh{\mra}\wh{\mrphi}$ (where $\wh{\mrphi}=\wh{\theta}\wh{\phi}$) holds. This implies that
\be\label{mrphi:mom}[(\otimes^j D)\otimes \wh{\mrphi}(\dm^{\tp}\cdot)](0)=[(\otimes^jD)\otimes(\wh{\mra}\wh{\mrphi})](0),\qquad j\in\N.\ee
Rearranging \er{mrphi:mom} yields
\be\label{mrphi:mom:1}\left[(\otimes^j D)\otimes \left(\wh{\mrphi}(\dm^{\tp}\cdot)^{\tp}-\wh{\mrphi}^{\tp}\wh{\mra}(0)^{\tp}\right)\right](0)=\left[(\otimes^jD)\otimes\left(\wh{\mrphi}^{\tp}(\wh{\mra}^{\tp}-\wh{\mra}(0)^{\tp})\right)\right](0),\ee
for all $ j\in\N$. Note that the right hand side of \er{mrphi:mom:1} only involves $\partial^\mu\wh{\mrphi}(0)$ with $|\mu|<j$. Furthermore, direct calculation yields
\be\label{mrphi:mom:2}\begin{aligned}&\left[(\otimes^j D)\otimes \left(\wh{\mrphi}(\dm^{\tp}\cdot)^{\tp}-\wh{\mrphi}^{\tp}\wh{\mra}(0)^{\tp}\right)\right](0)\\
	=&\left([(\otimes^jD)\otimes\wh{\mrphi}^{\tp}](0)\right)[(\otimes^j\dm^{\tp})\otimes I_r-(\otimes^jI_d)\otimes \wh{\mra}(0)],\end{aligned}\ee
for all $j\in\N$. Now by the condition in item (iii) on $\mra$, the matrix $[(\otimes^j\dm^{\tp})\otimes I_r-(\otimes^jI_d)\otimes \wh{\mra}(0)]$ is invertible for $j=1,\dots,m-1$. Moreover, it follows from \er{mrphi:mom:1} and \er{mrphi:mom:2} that up to a multiplicative constant, the values $\partial^\mu\wh{\mrphi}(0), |\mu|<m$ are uniquely determined via $\wh{\mrphi}(0)=\wh{\mra}(0)\wh{\mrphi}(0)$ and
$$\begin{aligned}&[(\otimes^jD)\otimes\wh{\mrphi}^{\tp}](0)\\
=&\left(\left[(\otimes^jD)\otimes\left(\wh{\mrphi}^{\tp}(\wh{\mra}^{\tp}-\wh{\mra}(0)^{\tp})\right)\right](0)\right)[(\otimes^j\dm^{\tp})\otimes I_r-(\otimes^jI_d)\otimes \wh{\mra}(0)]^{-1},\end{aligned}$$
for all $j=1,\dots,m-1.$

Consequently, by the above analysis and using \er{sr:mra} and \er{ref:mrta}, we conclude that \er{mrt:dual} holds for some $\tilde{C}\in\C\setminus\{0\}$, with $d\in\dlp{0}$ being a non-zero scalar multiple of $\mrd$.\\

On the other hand, the condition on $\mrta$ in item (iii) and item (v) together yield
\be\label{mrtphi:q}\wh{\mrtphi}(\xi)=K\wh{q}(\xi)\ol{\wh{\Vgu_{\dn}(\xi)}}^{\tp}+\bo(\|\xi\|^{\tilde{m}}),\qquad \xi\to 0,\ee
for some non-zero constant $K$. As item (ii) holds, then in particular item (3) of Theorem~\ref{thm:dfrt} holds. Then $\vmo(\tpsi)=\tilde{m}$ and \er{mrtphi:q} imply that $\wh{\mrtb}(\xi)\ol{\wh{\Vgu_{\dn}}(\xi)}^{\tp}=\bo(\|\xi\|^{\tilde{m}})$ as $\xi \to 0$. Now right multiplying $\wh{q}\ol{\wh{\Vgu_{\dn}}}^{\tp}$ to both sides of \er{dffrt} yields
\be\label{ref:mra}\ol{\wh{q}(\dm^{\tp}\xi)}\wh{\Vgu_{\dn}}(\dm^{\tp}\xi)\wh{\mra}(\xi)=\ol{\wh{q}(\xi)}\wh{\Vgu_{\dn}}(\xi)+\bo(\|\xi\|^{\tilde{m}}),\quad \xi\to 0.\ee
Since $\mra$ has $\tilde{m}$ sum rules with a matching filter $\wh{\mrvgu}:=\wh{\vgu}\wh{\theta}^{-1}$, we have $\wh{\mrvgu}(\dm^{\tp}\xi)\wh{\mra}(\xi)=\wh{\mrvgu}(\xi)+\bo(\|\xi\|^{\tilde{m}})$ as $\xi\to 0$. Moreover, $\mrta$ satisfies the refinement equation $\wh{\mrtphi}(\dm^{\tp}\cdot)=\wh{\mrta}\wh{\mrtphi}$. By the condition in item (iii) on $\mra$, we conclude from \er{mrtphi:q} and \er{ref:mra} that \er{mr:dual} must hold for some $C\in\C\setminus\{0\}$, with $c\in\dlp{0}$ being a non-zero scalar multiple of $q$.\\

Finally, by left multiplying $\ol{\wh{\mrphi}}^{\tp}$ and right multiplying $\wh{\mrtphi}$ (where $\wh{\mrtphi}:=\wh{\tilde{\theta}}\wh{\tphi}$) to \er{dffrt}, we have
$$\ol{\wh{\mrphi}(\dm^{\tp}\xi)}^{\tp}\wh{\mrtphi}(\dm^{\tp}\xi)=\ol{\wh{\mrphi}(\xi)}^{\tp}\wh{\mrtphi}(\xi)+\bo(\|\xi\|^{\tilde{m}+m}),\quad\xi\to 0.$$
By applying the same argument as in the proof of Lemma~\ref{lem:vguphi}, \er{mrphi:mrtphi} follows from the above identity. The proof is now complete.\qed\ep

\section{Summary}\label{sec:sum}

In this paper, we studied compactly supported multivariate OEP-based dual multiframelets with high order vanishing moments, and with a compact and banalced associated discrete multiframelet transform. We proved the main result Theorem~\ref{thm:dfrt} on the existence of such OEP-based dual multiframelets, with a constructive proof relying on a recently developed normal form of a matrix-valued filter. Furthermore, we provided structural analysis on compactly supported balanced OEP-based dual multiframelets.

Our investigation on OEP-based dual multiframelets focused on theoretical analysis. It is of practical interest to develop efficient algorithms to construct balanced dual multiframelets. However, this is well known that constructing multivariate framelets and wavelets are not easy in general. Moreover, the extremely strong conditions that both $\theta$ and $\tilde{\theta}$ must be strongly invertible makes the problem even harder. To achieve the strong invertibility on both filters $\theta$ and $\tilde{\theta}$, quite often it is unavoidable for them to have large supports, which is the main difficulty for us to perform construction in applications. Whether or not we can make the supports of both $\theta$ and $\tilde{\theta}$ as small as possible without sacrificing other desired properties is unknown. This could be a future research topic.


	%
	%

	

\end{document}